\documentclass[onefignum,onetabnum]{siamonline250211}

\usepackage{graphicx,epstopdf} 
\usepackage[caption=false]{subfig} 
\usepackage[ruled,vlined,algo2e]{algorithm2e}
\usepackage{mathtools}
\usepackage{amssymb}
\usepackage{wrapfig}
\usepackage{tikz}
\usepackage{booktabs}
\usetikzlibrary{trees, arrows.meta}

\newcommand {\Hi}{\mathcal{H}} 
\newcommand {\X}{\mathcal{X}} 
\newcommand {\R}{\mathbb{R}} 
\newcommand {\N}{\mathbb{N}} 
\DeclareMathOperator*{\argmax}{arg\,max}

\newtheorem{thm}{Theorem}[section]
\newtheorem{lem}[thm]{Lemma}
\newtheorem{cor}[thm]{Corollary}

\newtheorem{defn}[thm]{Definition}
\newtheorem{ass}[thm]{Assumption}

\usepackage{lipsum}
\usepackage{amsfonts}
\usepackage{graphicx}
\usepackage{epstopdf}
\usepackage{algorithmic}
\ifpdf
  \DeclareGraphicsExtensions{.eps,.pdf,.png,.jpg}
\else
  \DeclareGraphicsExtensions{.eps}
\fi

\usepackage{enumitem}
\setlist[enumerate]{leftmargin=.5in}
\setlist[itemize]{leftmargin=.5in}

\newsiamremark{remark}{Remark}
\newsiamremark{hypothesis}{Hypothesis}
\crefname{hypothesis}{Hypothesis}{Hypotheses}
\newsiamthm{claim}{Claim}
\newsiamremark{fact}{Fact}
\crefname{fact}{Fact}{Facts}

\headers{BO for inhomogeneous functions}{Olivier  Dondjio and Toni Karvonen}

\title{Bayesian Optimization for Inhomogeneous Functions via Locally Adaptive Kernels\thanks{
\funding{This work was supported by the Research Council of Finland grants 359183 (``Flagship of Advanced Mathematics for Sensing, Imaging and Modelling'') and 368086 (``Inference and approximation under misspecification'') as well as the Pilot for Doctoral Programmes by the Finnish Ministry of Education and Culture (pilot project ``Mathematics of sensing, imaging and modelling'').
}}}

\author{Olivier Dondjio\thanks{School of Engineering Sciences, Lappeenranta--Lahti University of Technology LUT, Finland
(\email{olivier.dondjio@lut.fi} \& \email{toni.karvonen@lut.fi}).}
\and Toni Karvonen\footnotemark[2]
}

\usepackage{amsopn}

\ifpdf
\hypersetup{
  pdftitle={Bayesian Optimization for Inhomogeneous Functions via Locally Adaptive Kernels},
  pdfauthor={Olivier Dondjio and Toni Karvonen}
}
\fi

\begin{document}

\maketitle

\begin{abstract}

Bayesian optimization (BO) is a powerful and general framework for sample-efficient global optimization.
Most BO methods rely on stationary Gaussian processes, implicitly assuming that an objective black-box function exhibits homogeneous smoothness across the input space. However, many complex real-world functions have highly inhomogeneous regularity, exhibiting both smooth and rough behavior. This fundamental mismatch between the model and the black-box function severely degrades both sample efficiency and surrogate accuracy of BO. To overcome this limitation, we introduce the Inhomogeneous Bayesian optimization framework. We develop a partitioned modeling strategy that dynamically subdivides the search space into subdomains, capturing spatially varying smoothness through locally adaptive kernels. To dynamically adjust the spatial resolution of the surrogate model, this approach utilizes a recursive splitting criterion. 
Furthermore, we develop a sampling allocation theory, establishing that driving the regional acquisition suprema below a common target level is achieved by polynomially denser sampling in rougher subdomains.
Numerical experiments across homogeneous functions, synthetic inhomogeneous functions, and a real-world 8-dimensional LASSO hyperparameter tuning task demonstrate that our approach improves optimization efficiency on synthetic inhomogeneous functions, while performing comparably to Classic BO on homogeneous functions and on the LASSO task.
This demonstrates the robustness of the adaptive partitioning mechanism.

\end{abstract}



\section{Introduction}

Bayesian optimization (BO) is a widely used approach for optimizing black-box functions that are expensive to evaluate, noisy, or analytically intractable. Since the early contributions of Mockus~\cite{mockus1998application} and the introduction of the efficient global optimization framework by Jones et al.~\cite{jones1998efficient}, BO has demonstrated strong performance in domains such as hyperparameter tuning of machine learning models~\cite{snoek2012practical}, robotics control~\cite{calandra2016bayesian}, materials discovery~\cite{gomez2018automatic}, and automated design~\cite{shahriari2015taking}. The standard BO framework combines a probabilistic surrogate, often a Gaussian process (GP), with an acquisition function that guides sampling by balancing exploration and exploitation.  

GPs provide a flexible prior distribution over functions with closed-form inference and uncertainty quantification~\cite{williams2006gaussian}. In the conventional BO setup, a GP is defined by a mean function (commonly assumed zero) and a covariance kernel encoding assumptions about smoothness and correlation between function values at different input locations. Stationary kernels such as the squared exponential and the Mat\'{e}rn family dominate the literature due to their simplicity and interpretability.

\begin{figure}[t] 
  \centering
  \vspace{-10pt} 
  
  \subfloat[1D Inhomogeneous function]{\label{f_inho_1}
    \includegraphics[width=0.48\textwidth]{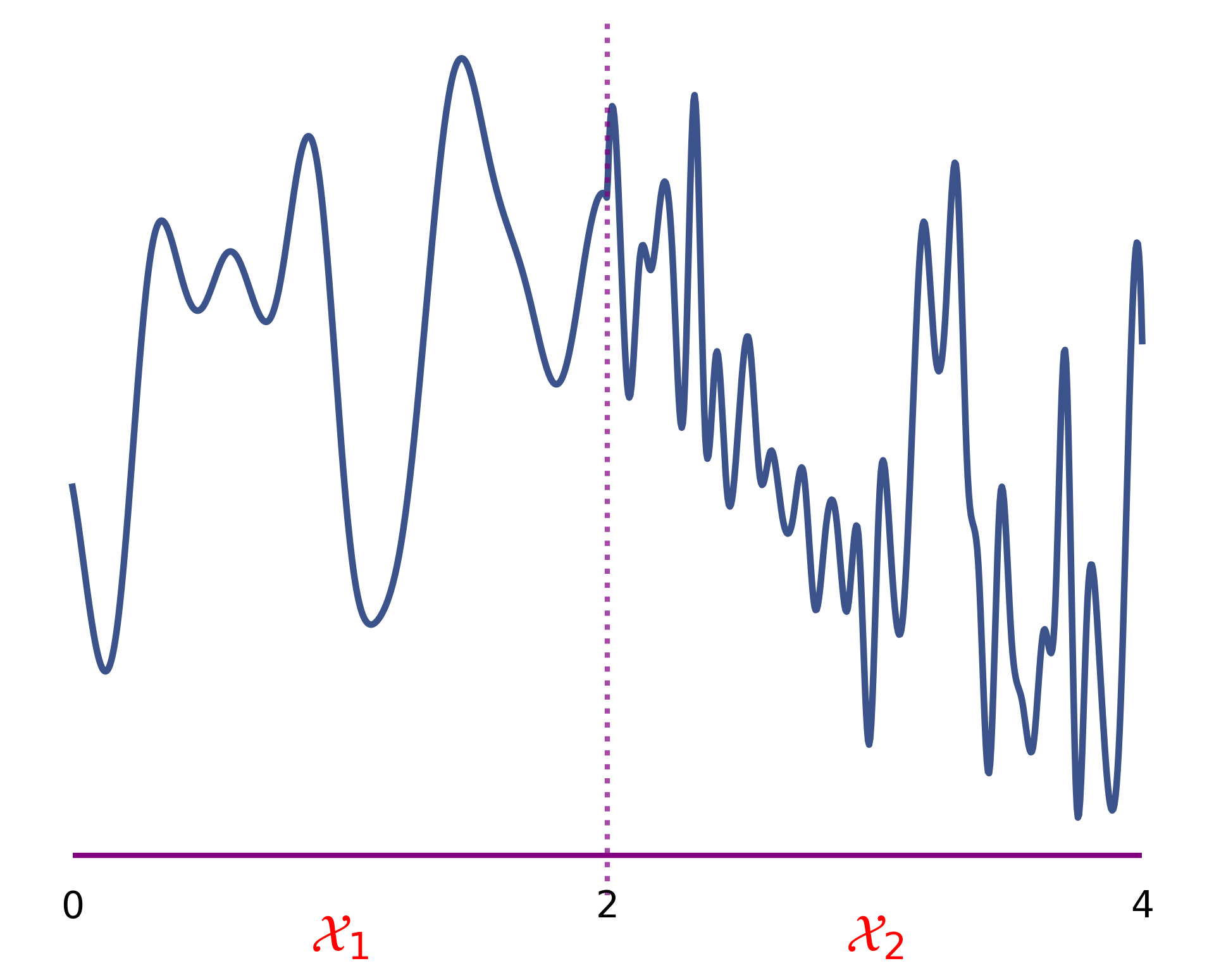}}
  \hspace{0.01\textwidth} 
  \subfloat[2D Inhomogeneous function]{\label{f_inho_2}
    \includegraphics[width=0.48\textwidth]{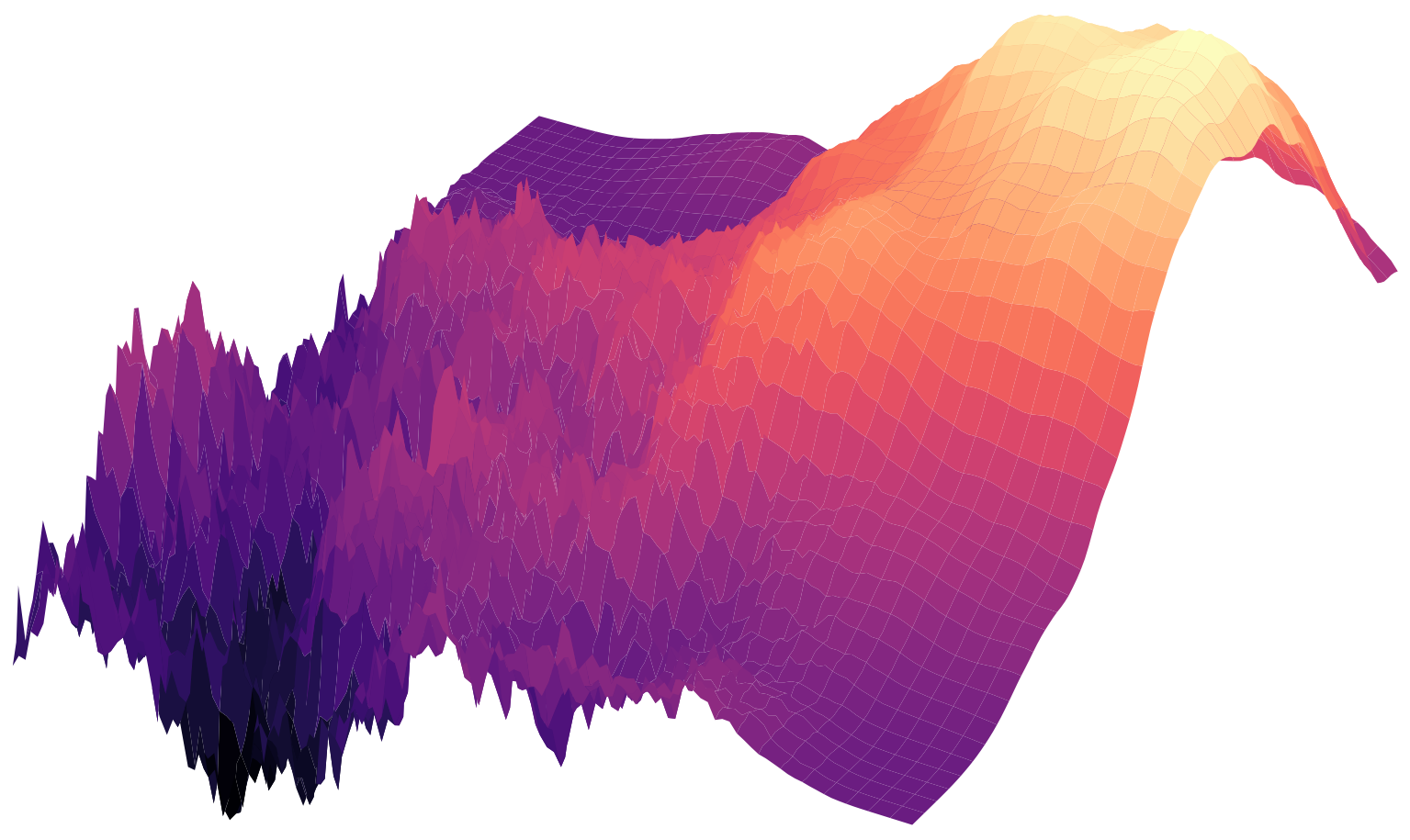}}

  \caption{Illustration of inhomogeneous behavior. In (a), the function shifts from a smooth regime in $ \X_1$ to high-frequency oscillations in $ \X_2$. In (b), whether or not the function is differentiable depends on the spatial location.}  
  \label{f_inho}
\end{figure}

However, stationary kernels assume that the statistical properties of the function are invariant under translation in the input space. In other words, using them amounts to assuming that the objective function has the same smoothness everywhere on its domain. This assumption is violated when the objective function exhibits heterogeneous smoothness: as illustrated in Figure~\ref{f_inho}, the function may vary smoothly in certain subdomains while exhibiting sharp transitions, high-frequency oscillations, or local irregularities in others.
This article develops new BO methodology for inhomogeneous functions such as those in Figure~\ref{f_inho}.

\subsection{Related Work}
 
To address such settings, various locally adaptive extensions of GPs have been proposed.
These can generally be categorized into domain partitioning, input transformation, and localized non-stationary modeling. In the domain partitioning category, Gramacy and Lee~\cite{gramacy2008bayesian} developed treed GPs that partition the search space via a data-driven tree structure to fit independent local GPs with region-specific hyperparameters, while more recently Eriksson et al.~\cite{eriksson2019scalable} proposed TuRBO, which maintains a collection of simultaneous local trust regions with independent GP models, avoiding the need to model global heterogeneity entirely. 
In the localized-modeling category, Paciorek and Schervish~\cite{paciorek2003nonstationary} introduced a general construction for
explicitly non-stationary covariance functions, subsequently applied by Plagemann et al.~\cite{plagemann2008nonstationary}, whose non-stationary GPs estimate kernel hyperparameters from local point estimates of smoothness across a single, unpartitioned domain; related non-stationary constructions include \cite{higdon1999non,schmidt2003bayesian}. 
Finally, input-transformation approaches instead reshape the input space before applying a standard stationary kernel. Snoek et al.~\cite{snoek2014input} propose a fixed, monotonic per-dimension warp via the beta distribution, while Calandra et al.~\cite{calandra2016manifold} generalize this to a jointly learned, non-monotonic feature-space mapping (manifold GPs), capable of capturing discontinuities that a monotonic warp cannot represent.

Recent studies reinforce the importance of local adaptation in GP-based modeling. Naslidnyk et al.~\cite{naslidnyk2025comparing} and Karvonen and Bachoc~\cite{KarvonenBachoc2025} compared scale parameter estimators for GP interpolation under Brownian motion and periodic Matérn priors, showing that leave-one-out cross-validation provides better-calibrated uncertainty estimates than maximum likelihood estimation across a broad class of functions, highlighting the role of local parameter selection. Liautaud et al.~\cite{pmlr-v272-liautaud25a} introduced a minimax adaptive online regression algorithm over Besov spaces, capable of dynamically adjusting prediction resolution according to spatially inhomogeneous smoothness. Their results demonstrate that locally adaptive strategies substantially improve performance in heterogeneous environments, aligning conceptually with the motivation for adaptive kernels in BO.  

A separate line of work concerns convergence rates for Bayesian optimization rather than the surrogate model itself. Bull~\cite{bull2011convergence} established minimax rates for expected improvement over Sobolev balls, and
Vazquez and Petit~\cite{vazquez2026simple} have recently sharpened the upper bound for the fixed-prior policy. 
However, these analyses assume that the objective lies in the Reproducing Kernel Hilbert Space (RKHS) of the kernel. The framework developed here concerns the complementary case, in which the local regularity of the objective differs from the regularity the kernel assumes and the scale parameter is estimated locally from regional observations.

\subsection{Methodology and Contributions}

To address spatial heterogeneity without sacrificing the analytical tractability of standard stationary kernels, Section~\ref{sec:IBO} introduces the \emph{Inhomogeneous Bayesian Optimization} (IBO) framework. 
The core idea is to subdivide the input domain $\X$ into subdomains $\X_i$ and use a locally scaled kernel of the form
\[
k_{\sigma_i}(x, x') = \sigma_i^2 k(x, x') \quad \text{ for } \quad x, x' \in \X_i,
\]
where the local scale parameter $\sigma_i$ controls the regional prior signal variance and predictive uncertainty across subdomains. By estimating $\sigma_i$ locally within each partition using local observations, the model expands uncertainty bounds in rougher, high-variance subdomains while contracting them in smoother, low-variance zones.
Importantly, local scale parameters do not affect the posterior mean, which remains continuous across $\X$ whenever $k$ is continuous.

Compared to the locally adaptive methods discussed above, the proposed IBO framework distinguishes itself in two primary ways. First, unlike global non-stationary covariance functions that require optimizing complex, spatially varying hyperparameters, IBO relies on a partition-based approach that maintains the robust, closed-form properties of stationary kernels locally. Second, while conceptually related to treed GPs~\cite{gramacy2008bayesian}, our framework is explicitly designed to streamline the BO loop. It replaces heavy tree-search inference with a lightweight recursive splitting criterion based on local scales. By directly integrating these local scale estimates into localized acquisition functions, IBO provides a mathematically grounded framework for sampling allocation that standard global or treed methods lack.

To formalize this allocation strategy, Section~\ref{sec:SATIBO} establishes a rigorous theoretical foundation for sequential sampling in inhomogeneous BO. We prove (Theorems~\ref{th:thm0} and~\ref{cor:cor42}) that a uniform upper bound on the acquisition function across all subdomains is achieved if local sample sizes obey polynomial scaling relations, with exponents determined by the ratios of local decay rates.
Polynomially more samples in subdomains of lower smoothness suffice to achieve this.

Finally, in Section~\ref{sec:NIL} we validate the practical efficacy of the framework by evaluating IBO across standard homogeneous benchmark functions, synthetic inhomogeneous functions, and a real-world $8$-dimensional sparse regression hyperparameter tuning problem (LassoBench). The numerical experiments demonstrate that IBO outperforms Classic BO on inhomogeneous objective landscapes by accelerating convergence. 
Crucially, when applied to standard homogeneous functions, the criterion divides those objectives whose local scale varies across the domain, and IBO still achieves performance comparable to standard global methods.

\section{Bayesian Optimization} 
\label{sec:BO}

Let \( f \colon  \X \to \R \) be a black-box objective function on \mbox{$\X \subset \R^d$} that we aim to maximize.
We assume that evaluations of \( f \) are expensive and that we lack access to analytic gradients.
BO provides a principled framework for finding the global optimizer of such functions with as few evaluations as possible~\cite{brochu2010tutorial, shahriari2015taking}.  Rather than relying on exhaustive or unguided search strategies, BO constructs a probabilistic surrogate model to approximate the objective function and an acquisition function that guides the search for an optimum. The surrogate model, most commonly a GP model, provides both predictions and uncertainty estimates of \( f \)~\cite{williams2006gaussian}. The acquisition function then leverages this information, specifically the posterior mean and variance of the GP, to determine the most promising point  for the next evaluation by balancing exploration of uncertain regions with exploitation of areas that are  likely to yield improvement~\cite{brochu2010tutorial,jones1998efficient}.

\subsection{Gaussian Process Regression}\label{SubSec:GPR}

A Gaussian process $f \sim \mathcal{GP}(m, k)$ is a distribution over functions such that any finite collection of function values follows a multivariate Gaussian distribution~\cite{williams2006gaussian}. It is fully specified by a mean function $m \colon \X \to \mathbb{R}$ and a covariance kernel $k \colon  \X \times  \X \to \mathbb{R}$. The mean function represents the expected value of the function at any given point, $m(x) = \mathbb{E}[f(x)]$, while the covariance kernel captures the correlation between function values at different points, $k(x, x') = \mathbb{E}[\big(f(x) - m(x)\big)\big(f(x') - m(x')\big)]$.

As is standard, we assume a zero prior mean function (i.e., $m = 0$). Then our prior beliefs about the function's smoothness and correlation structure are encoded entirely within the covariance kernel. One of the most commonly used classes of kernels in BO is the Mat\'{e}rn class of kernels. Its popularity stems from its ability to model functions with varying degrees of differentiability, offering a more flexible and realistic alternative to the infinitely smooth squared exponential kernel~\cite{KarvonenOates2026}.
For $x, x' \in \R^d$, the Mat\'{e}rn covariance kernel is defined as
\begin{align}
    k_{\text{Mat\'{e}rn}}(x, x') = \frac{2^{1-\nu}}{\Gamma(\nu)} \left( \sqrt{2\nu} \frac{\|x - x'\|}{\ell} \right)^\nu K_\nu \left( \sqrt{2\nu} \frac{\|x - x'\|}{\ell} \right),\label{eq:matern}
\end{align}
where $K_\nu$ is the modified Bessel function of the second kind. The properties of the sample paths generated by this kernel are governed by two parameters:
\begin{itemize}
    \item \textbf{Smoothness parameter $\nu$}: This is the defining feature of the Mat\'{e}rn family. It directly controls the mean-square differentiability of the GP; a process is $s$-times mean-square differentiable if and only if $\nu > s$~\cite[Section~2.7]{stein1999interpolation}. Common choices include $\nu = 3/2$ (once differentiable), $\nu = 5/2$ (twice differentiable) and the limiting case $\nu \to \infty$, which recovers the infinitely differentiable squared exponential kernel.
    \item \textbf{Length-scale $\ell$}: This parameter controls how quickly correlations between function values decay in the input space. A large $\ell$ implies that the function values change slowly across the input space, while a small $\ell$ allows for more rapid oscillations.
\end{itemize}
Any kernel, including Matérns, can be equipped with an explicit scale parameter that controls the overall magnitude of the variance:
\[
k_\sigma(x,x') = \sigma^2 k(x,x'),
\]
where $k$ is a base kernel and $\sigma^2$ determines the prior variance of the process. 
A large scale parameter yields a GP prior $f \sim \mathcal{GP}(0, k_\sigma)$ with large variance.
As explored by Karvonen et al.~\cite{karvonen2020maximum} and Naslidnyk et al.~\cite{naslidnyk2025comparing}, one of the advantages of maintaining an explicit scale parameter lies in its ability to compensate for model misspecification.

Suppose we have evaluated the objective function at a set of points 
$X = \{x_1, \dots, x_{N}\} \subset \X$.
We denote the corresponding function values as
\[
f(X) = \big(f(x_1), \dots, f(x_{N})\big)^\top.
\]
Under the GP prior $f \sim \mathcal{GP}(0, k_\sigma)$, these observations follow a multivariate normal distribution:
$$f(X) \sim \mathcal{N}\big(\mathbf{0}, \sigma^2 K(X,X)\big),$$
where $K(X,X) \in \mathbb{R}^{N \times N}$ is the kernel matrix with entries 
\(
[K(X,X)]_{ij} = k(x_i, x_j),
\)
which encodes the correlations among the observed inputs.
Given the dataset $\mathcal{D}_{N} = \{(x_i, f(x_i))\}_{i=1}^{N}$, the GP posterior distribution at a new point $x$ is obtained by conditioning the joint Gaussian prior on the observed data $f(X)$. 
We denote the posterior mean at $x \in \mathcal{X}$ conditioned on $\mathcal{D}_N$ by $m(x ; \mathcal{D}_N)$. To simplify notation, we write $m_N(x) = m(x ; \mathcal{D}_N)$ whenever the dataset is clear from the context.
The posterior mean is
\begin{align}
    m_N(x) = k(x,X)^\top K(X,X)^{-1} f(X)\label{eq:mn}.
\end{align}
Here $k(x,X) = (k(x,x_1), \dots, k(x,x_{N}) )^\top$ is the vector of kernel evaluations between a new point $x$ and the training inputs.
The posterior mean does not depend on the scale parameter~$\sigma$.
The posterior covariance is
\begin{align*}
    k_{\sigma,N}(x,y) =\sigma^2 k_{N}(x,y),
\end{align*}
where $k_N(x,y) = k_N(x,y ; \mathcal{D}_N)$ is the posterior covariance corresponding to the unscaled base kernel $k$:
\begin{align*}
k_{N}(x,y) = k(x,y) - k(x,X)^\top K(X,X)^{-1} k(X,y) .
\end{align*}
We use the shorthands $k_{\sigma,N}(x)=k_{\sigma,N}(x,x)$ and $k_{N}(x) = k_{N}(x,x)$ for the posterior variance.

\begin{figure}[t]
\centering
\includegraphics[width=\textwidth]{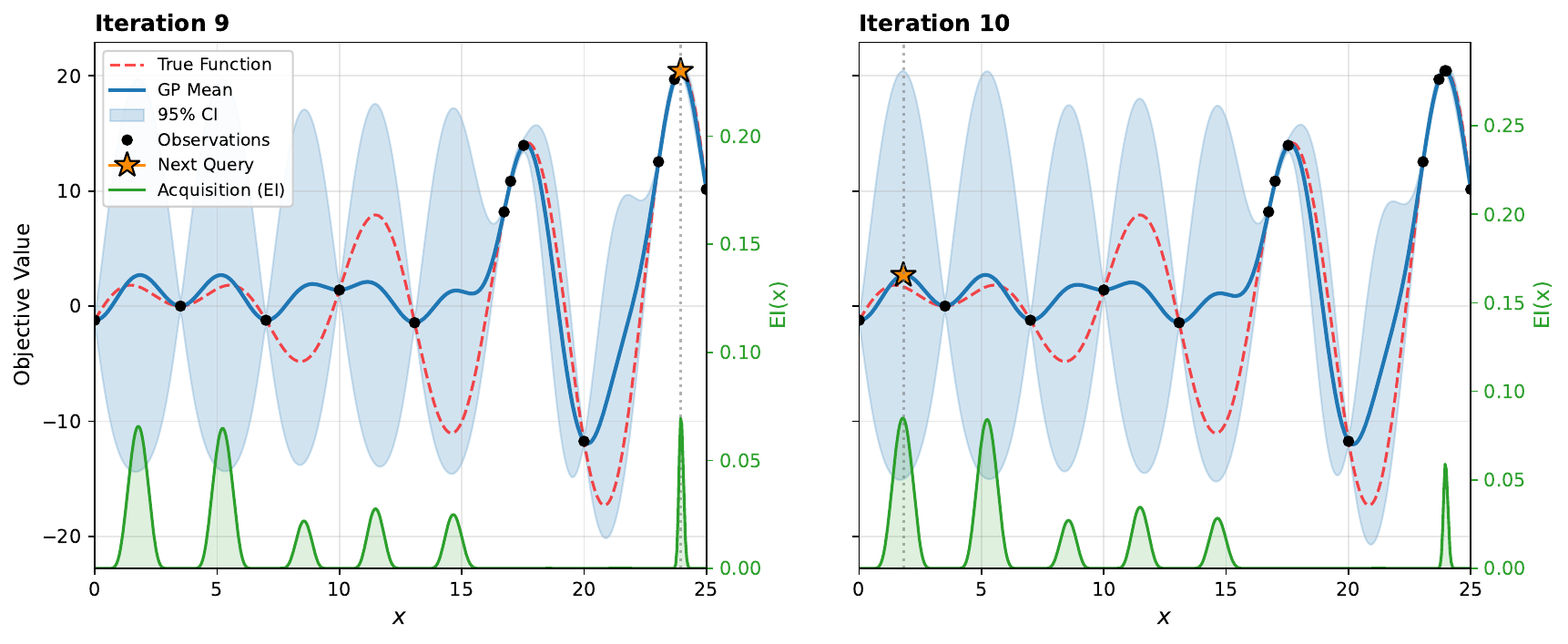}
\caption{Classic BO on a univariate function. The shaded region represents the 95\% credible interval, and the green curve indicates the EI.}
\label{fig:BO}
\end{figure}
The resulting posterior distribution at any point $x$ is Gaussian with mean $m_N(x)$ and variance $\sigma^2 k_N(x)$.
The posterior mean serves as the surrogate model's prediction, while $\sigma^2 k_N$ quantifies the remaining uncertainty after observing the data. At this stage, the scale factor $\sigma^2$ acts as a global multiplier, uniformly scaling the overall magnitude of this uncertainty across the entire domain, while the base variance $k_N$ governs how uncertainty naturally increases as we move further away from the observed data points. This behavior is illustrated in Figure~\ref{fig:BO}, where the solid line represents the posterior mean $m_N(x)$ and the shaded region the $95\%$ credible interval $[m_N(x) - 1.96 \sigma \sqrt{k_N(x)}, m_N(x) + 1.96  \sigma \sqrt{k_N(x)}]$ derived from the posterior variance.
GP regression is well suited for BO, as it provides both predictions and calibrated uncertainty estimates in closed form, which are crucial for defining acquisition functions.

\subsection{Acquisition Function}\label{subsec:AC}

The acquisition function $a (\cdot; \mathcal{D}_N) \colon \X \to \mathbb{R}$ is a key component of BO, as it determines where the objective function should be evaluated next. By leveraging the surrogate model, the acquisition function quantifies the utility of sampling at any point $x$ by balancing exploration of uncertain regions and exploitation with optimal predicted values. The next evaluation point $x_{N+1}$ is selected by maximizing the acquisition function: 
$$x_{N+1} =  \underset{x\in \X}{\argmax}~a(x; \mathcal{D}_{N}).$$
Several acquisition functions have been proposed.
Among the most popular are Expected Improvement (EI)~\cite{jones1998efficient}, the Probability of Improvement (PI)~\cite{kushner1964new}, the Lower Confidence Bound (LCB)~\cite{srinivas2009gaussian}, 
and Thompson Sampling (TS)~\cite{thompson1933likelihood}, which is a stochastic acquisition strategy. 

The PI focuses on selecting points that maximize the 
probability of outperforming the current best observation, which can lead to  over-exploitation. In contrast, the LCB explicitly balances exploration and exploitation by combining the posterior mean and variance through a  tunable parameter. TS, on the other hand, adopts a stochastic approach by drawing samples from the 
posterior distribution and selecting points that are optimal with respect to these  samples, thereby naturally incorporating exploration. Each of these methods offers a different trade-off between exploration and exploitation, 
which is a central challenge in BO.

In this work, we focus our theoretical analysis on the $\mathbf{EI}$ acquisition function due to its analytical tractability, while our numerical experiments utilize all the aforementioned acquisition functions.
 The $\mathbf{EI}$ criterion is one of the most widely used acquisition functions in BO~\cite{brochu2010tutorial,jones1998efficient}. Its purpose is to guide the search for the global optimum by identifying points in the input space that are likely to yield the greatest progress relative to the best value observed so far.
At iteration $N$, suppose the best function value observed is
\[
f(x^*) = \max \{ f(x_1), \dots, f(x_{N}) \}.
\]
For a candidate point \(x\), the improvement over the current best is defined as
\[
\mathbf{I}(x;\mathcal{D}_{N}) = \max(f(x) - f(x^*), 0).
\]
The $\mathbf{EI}$ at point $x$ is then defined as the expectation of this improvement under the posterior distribution of the surrogate model:
\begin{align*}
    \mathbf{EI}(x;\mathcal{D}_{N}) = \mathbb{E}[\mathbf{I}(x;\mathcal{D}_{N})].
\end{align*}
By evaluating this expectation using the GP posterior mean $m_N(x)$ and variance $k_{\sigma, N}(x)$, the $\mathbf{EI}$ criterion measures the expected amount by which evaluating $x$ would improve over $f(x^*)$. This formulation naturally balances \emph{exploitation}, which favors regions where the predicted mean is higher than \(f(x^*)\), and \emph{exploration}, which favors regions with large uncertainty, even if the predicted mean appears suboptimal.

An analytic expression for $\mathbf{EI}$ is obtained by integrating the improvement function against the Gaussian probability density of the surrogate model at $x$~\cite{jones1998efficient}:
\begin{align*}
\mathbf{EI}(x;\mathcal{D}_{N}) =
\begin{cases}
(m_N(x)-f(x^*))\,\Phi\big(z(x)\big)
+ \sigma\sqrt{k_{N}(x)}\,\phi \big(z(x)\big) & \text{ if } \quad k_{N}(x)>0,\\[1ex] 
0 & \text{ if } \quad k_{N}(x)=0 .
\end{cases}
\end{align*}
Here,
\[
z(x) = \frac{m_N(x) - f(x^*)}{\sigma\sqrt{k_{N}(x)}}.
\]
The functions $\Phi$ and $\phi$ are the cumulative distribution function and the probability density function of the standard normal distribution, respectively.
The formula consists of two additive terms. The first is the \emph{predicted improvement term}, 
\((m_N(x) - f(x^*) )\, \Phi(z(x))\), which accounts for the improvement expected from the mean prediction. 
The second is the \emph{uncertainty term}, 
\(\sigma \sqrt{k_{N}(x)} \, \phi(z(x))\), which captures the potential gain due to variance and thereby encourages exploration.
The next evaluation point is determined by maximizing $\mathbf{EI}$ (see Figure~\ref{fig:BO}):
\begin{equation}
    x_{N+1} = \underset{x\in \X}{\argmax}  \, \mathbf{EI}(x;\mathcal{D}_{N}).\label{eq:Nextpoint}
\end{equation}
In practice, this maximization is often performed using a gradient-based 
optimization method, since $\mathbf{EI}$ is differentiable with respect to $x$ wherever $k_N(x) > 0$. 
This allows efficient navigation of the acquisition surface and ensures rapid identification of the next query point~\cite{jones1998efficient}.

\section{Inhomogeneous Bayesian Optimization} \label{sec:IBO}

Standard BO methods based on commonly used stationary kernels typically assume \textit{homogeneity} in the objective function's smoothness, adopting a global kernel with a single smoothness or length-scale parameter.
However, in many real-world applications, the function $f$ that is being optimized exhibits heterogeneous regularity across the input space $\X$. As illustrated in Figure~\ref{f_inho}, such functions may be characterized by regions of slow variation contrasted with areas of rapid fluctuation or sharp transitions. Applying a global surrogate model to these \textit{inhomogeneous} functions often leads to inefficiencies: the model over-smooths high-frequency features or under-smooths stable regions.
To address this, we introduce the Inhomogeneous Bayesian optimization framework. The core principle of IBO is to allow the statistical properties of the surrogate model to adapt locally in a way that retains the ease-of-use of GPs.
This is particularly relevant in fields such as fluid dynamics and cavitation studies~\cite{bulanov2023cavitation,lu2025spectral,sammuller2023neural}.
We model heterogeneity by decomposing the domain $ \X$ into $n$ subdomains:
\begin{align*}
    \X &= \bigcup_{i=1}^{n}\X_i ,\label{eq:partition_sum}
\end{align*}
where the subdomains have pairwise disjoint interiors.
To ensure the optimization is well-defined, each subdomain $\X_i$ is closed.
While adjacent subdomains can overlap on their boundaries, any input point $x$ falling exactly on a shared boundary is uniquely assigned to the subdomain with the smallest index $i$. 
This partitioning allows us to define the surrogate model piecewise: for instance, in Figure~\ref{f_inho_1}, the domain is split into $\X_1$ and $\X_2$ based on the point at which the function's behavior changes.
The inputs in $\X_i$ are denoted $X_i = X \cap \X_i^* \subset X$, where $\X_i^* = \X_i \setminus \cup_{k=1}^{i-1} \X_k$, and the local dataset is $\mathcal{D}^i_{N} = \{(x_j, f(x_j))  : x_j \in X_i\} \subset \mathcal{D}_N$.
We write $N_i = |X_i|$ for the number of input points located in $\X_i$.
Within each partition $\X_i$, the IBO framework parameterizes the surrogate model with a local scale parameter $\sigma_i$, allowing the model to accurately reflect the regularity of $f$ in that region. This approach is conceptually related to treed GPs~\cite{gramacy2008bayesian}, but optimized for the sequential decision-making process of BO.

The primary objective behind this localized parameterization is to dynamically manage the sampling density across the domain. By assigning a unique scale parameter $\sigma_i$ to each subdomain $\X_i$, the model can distinguish between regions of high and low complexity. Our goal is to use these local estimates to inform the sampling allocation: we aim to place more evaluation points in subdomains where the function exhibits higher variability, while maintaining a sparser set of observations in subdomains in which $f$ is smoother. This ensures that the global search budget is not wasted on homogeneous areas and is instead concentrated where the local characteristics of $f$ demand higher spatial resolution.

\subsection{Single-kernel case}\label{sec41}

We begin with a baseline model where the kernel family $k$ is fixed globally, but the scale parameter $\sigma$ is allowed to vary across subdomains to capture the local properties of $f$. Specifically, the magnitude of this local scale parameter reflects the function's regularity: rougher regions are characterized by larger effective variability (a larger local scale), whereas smoother regions exhibit smaller variability. 
This is supported by the theoretical findings in~\cite{KarvonenBachoc2025, karvonen2020maximum}.
As established in~\eqref{eq:mn}, the scale parameter does not affect the posterior mean. 
The GP posterior distribution at a new point $x \in \X$ has mean $m_N(x)$ given in~\eqref{eq:mn}, while the posterior covariance in $\X_i$ is
\begin{align*}
    k_{\hat{\sigma}_i,N}(x,y) = \hat{\sigma}_i^2 k_{N}(x,y)  \quad \text{for } x,y\in \X_i,
\end{align*}
where $\hat{\sigma}_i$ is a local estimator of the scale parameter $\sigma$, computed using only the subset of training data that falls within the corresponding subdomain $ \X_i$. One way to estimate $\sigma$ is the Maximum Likelihood (ML) estimator~\cite[Section~2]{jones1998efficient} defined as follows:
\begin{align}\label{eq:sigma}
    \hat{\sigma}_{i}^2
    = \frac{f(X_i)^{\top} K(X_i,X_i)^{-1} f(X_i)}{N_i}.
\end{align}

We define the global acquisition function $a(x;\mathcal{D}_{N})$ piecewise across the subdomains. This allows the search strategy to account for regional differences in function regularity. For a partition $\{\X_i\}_{i=1}^n$ of $\X$, the global acquisition function is defined as
\begin{align}
a(x;\mathcal{D}_{N}) = 
\begin{cases} 
a_1(x;\mathcal{D}_{N}), & x \in \X_1, \\
\vdots \\
a_n(x;\mathcal{D}_{N}), & x \in \X_n,
\end{cases} \label{eq:EI}
\end{align}
where $a_i(x;\mathcal{D}_{N})$ denotes the acquisition function computed using the local scale and kernel parameters specific to region $\X_i$.

As a primary example, we consider the localized Expected Improvement.
For each region $\X_i$ and each $x$ with $k_N(x) > 0$, define the locally standardized score
\begin{align}
z_i(x) = \frac{m_{N}(x) - f(x^*)}
{\hat{\sigma}_i \sqrt{k_{N}(x)}}. \label{eq:zi}
\end{align}
The regional acquisition function is then
\begin{align}
\mathbf{EI}_i(x;\mathcal{D}_{N}) =
\begin{cases}
  \big(m_{N}(x) - f(x^*)\big)\,\Phi\!\left(z_i(x)\right)
  + \hat{\sigma}_i \sqrt{k_{N}(x)} \,\phi\!\left(z_i(x)\right), & \text{ if } \quad k_{N}(x) > 0, \\[1ex]
  0, & \text{ if } \quad k_{N}(x) = 0.
\end{cases} \label{eq:EIi}
\end{align}
By substituting $\mathbf{EI}_i(x;\mathcal{D}_{N})$ into the piecewise framework \eqref{eq:EI}, we obtain a global acquisition function $\mathbf{EI}(x;\mathcal{D}_{N})$.

\subsection{Multi-kernel case}

In the single-kernel configuration discussed in Section~\ref{sec41}, we assumed that all subdomains 
\(\{\X_i\}_{i=1}^n\) share the same kernel \(k\) and are differentiated only by region-specific scale parameters $\sigma_i$.
While this approach is effective when the objective function $f$ exhibits varying local magnitudes within a consistent functional class, many practical scenarios involve regions with fundamentally different structural properties. For instance, a function might appear nearly linear in one subdomain, oscillate rapidly in another, or contain sharp, localized peaks in a third.
Because a single kernel cannot effectively capture such diverse local structures, a more flexible approach is sometimes required.
To address these diverse behaviors, one could assign each region $ \X_i$ its own distinct kernel $k_i$ alongside its local scale parameter $\sigma_i$.
All other quantities introduced in the single-kernel setting, including the posterior covariance and the localized Expected Improvement criterion, would then be defined analogously using the local kernel \(k_i\) and the local dataset \(\mathcal{D}^i_{N}\).

A primary challenge in the multi-kernel setting is the interface problem. Because each local posterior mean is constructed using only local data $\mathcal{D}^i_{N}$ and a local kernel $k_i$, there is no inherent mechanism to align predictions over the boundaries of the subdomains.
For $d=1$, it is often possible to enforce continuity across subdomains through boundary constraints or shared boundary nodes. However, for $d > 1$, the increased dimensionality and the complexity of the boundaries, which are surfaces or hyperplanes, make enforcing global continuity significantly more difficult.

\subsection{Recursive Partitioning}\label{subsec:RPHAO}

To adapt the spatial resolution of the model to the local variability of the objective function, the input domain $ \X$ is recursively partitioned into a collection of leaf regions $\mathcal{L}=\{\mathcal{X}_1,\mathcal{X}_2,\dots,\mathcal{X}_n\}$ with disjoint interiors using a top-down binary splitting strategy evaluated at the median of the observations along each coordinate axis. In contrast to uniform decompositions or Voronoi-based partitioning schemes, this adaptive construction allocates finer subdivisions to highly heterogeneous areas while preserving larger regions in approximately stationary zones.

The informativeness of a candidate split is assessed through local scale.
For a region $\mathcal{X}_l$ containing $N_l$ observations, the local scale $\hat{\sigma}_l$
is estimated using an ML-based estimator~\eqref{eq:sigma}.

\begin{figure}[t]
\centering

\begin{minipage}[c]{0.54\textwidth} 
\begin{algorithm}[H]
\DontPrintSemicolon
\SetAlgoLined
\KwIn{Domain $\X$, Current observations $\mathcal{D}_N$, threshold $\tau$, max depth $D_{\mathrm{max}}$}
\KwOut{Partition leaves $\mathcal{L}$}
Initialize $\mathcal{L} = \{\X\}$, $D = 1$\;\BlankLine
\While{$D \le D_{\mathrm{max}}$}{\BlankLine
    $\text{split\_occurred} \leftarrow \text{False}$\;\BlankLine
    $\mathcal{C} \leftarrow \mathcal{L}$ \BlankLine
    \For{$\mathcal{R} \in \mathcal{C}$}{
        \For{each axis $a \in \{1, \dots, d\}$}{
            Split $\mathcal{R}$ at the median into $\mathcal{X}_{\text{left}, a}$ and $\mathcal{X}_{\text{right}, a}$\;
            \BlankLine
            Compute local scales $\hat{\sigma}_{\text{left},a}$ and $\hat{\sigma}_{\text{right},a}$\;
            \BlankLine
            Compute $\eta_a$\;
        }
        $a^* \leftarrow \argmax_a \eta_a$\;\BlankLine
        \If{$\eta_{a^*} > \tau$}{
            Remove $\mathcal{R}$ from $\mathcal{L}$\;
            \BlankLine
            Add $\mathcal{X}_{\text{left}, a^*}$ and $\mathcal{X}_{\text{right}, a^*}$ to $\mathcal{L}$\;
            $\text{split\_occurred} \leftarrow \text{True}$\;
        }
    }
    \If{\textbf{not} $\text{split\_occurred}$}{
        \textbf{break} 
    }
    $D \leftarrow D + 1$\;
}
\Return{$\mathcal{L}$}
\caption{Recursive Partitioning Strategy}
\label{alg:partitioning}
\end{algorithm}
\end{minipage}\hfill
\begin{minipage}[c]{0.44\textwidth} 
    
\centering
\resizebox{\textwidth}{!}{
\begin{tikzpicture}[
    level distance=1.3cm,
    level 1/.style={sibling distance=4.2cm},
    level 2/.style={sibling distance=2.2cm}, 
    every node/.style={
        rectangle,
        draw,
        rounded corners,
        text width=1.7cm,
        align=center,
        fill=gray!10,
        font=\scriptsize,
        inner sep=2pt
    },
    edge from parent/.style={draw, ->, thick}
]

\node {Root $\mathcal{X}_0= \X$ \\ $\hat{\sigma}_0$}
    child {
        node {Split $a^*$ \\ $\eta_{a^*} > \tau$}
        child {
            node {$\mathcal{X}_{\mathrm{left}}$ \\ $\hat{\sigma}_{\mathrm{left}}$}
            child {node[fill=green!20] {Leaf \\ $\eta_{a^*} \leq \tau$}}
            child {node {Split ...}}
        }
        child {
            node {$\mathcal{X}_{\mathrm{right}}$ \\ $\hat{\sigma}_{\mathrm{right}}$}
            child {node {Split ...}}
            child {node[fill=green!20] {Leaf \\ $D=D_{\mathrm{max}}$}}
        }
    };
\end{tikzpicture}
}
\caption{Hierarchical structure of the partition tree. Green nodes represent termination points where the region is sufficiently stationary or the maximum depth is reached.}
\label{fig:partition_tree}
\end{minipage}

\end{figure}

To determine the optimal partitioning direction, candidate splits are evaluated along each coordinate axis $a \in \{1,\dots,d\}$.
For each axis, the region is temporarily divided into two subregions, $\mathcal{X}_{\text{left},a}$ and $\mathcal{X}_{\text{right},a}$, and their corresponding local scales $\hat{\sigma}_{\text{left},a}$ and $\hat{\sigma}_{\text{right},a}$ are estimated.
To quantify the degree of heterogeneity along any candidate axis $a$, we compute the \emph{non-stationarity ratio} $\eta_a$, defined as
\begin{equation*}
\eta_a = \frac{\max(\hat{\sigma}_{\text{left},a}, \hat{\sigma}_{\text{right},a})}{\min(\hat{\sigma}_{\text{left},a}, \hat{\sigma}_{\text{right},a})}.
\end{equation*}
The optimal splitting axis $a^*$ is chosen as the one that maximizes this relative non-stationarity ratio, rather than merely the largest one-sided scale:
\begin{equation*}
a^* = \argmax_{a} \eta_a.
\end{equation*}
This ratio serves as a diagnostic tool: a high maximum non-stationarity ratio $\eta_{a^*}$ suggests that the parent region covers functional behaviors with significantly different magnitudes, justifying the division into independent local models. 

A split is performed only if this maximum ratio exceeds a predefined non-stationarity threshold $\tau$, set to $1.5$ in the experiments of Section~\ref{sec:NIL}.
This ensures that the same heterogeneity criterion governs both axis selection and split acceptance.
To prevent excessive fragmentation and limit the computational overhead of managing too many local models, the recursive partitioning is constrained by a second stopping condition: a maximum tree depth, denoted $D_{\mathrm{max}}$ and set to $3$ in Section~\ref{sec:NIL}. A candidate split is rejected and the region becomes a terminal leaf if  executing the split would cause the depth of that branch in the partition tree to exceed $D_{\mathrm{max}}$. Regions exhibiting strong relative variations in local scale are recursively subdivided, producing a dense local grid, whereas approximately stationary regions remain large, thereby preserving computational efficiency.

The overall recursive procedure used to construct the partition tree is summarized in Figure~\ref{fig:partition_tree} and Algorithm~\ref{alg:partitioning}. 
The resulting partition tree defines a collection of local Gaussian process models, each associated with a
subdomain $\X_i$ and its corresponding regional $\mathbf{EI}$ acquisition function $\mathbf{EI}_i(x;\mathcal{D}_{N})$. The local scale estimate
$\hat{\sigma}_i$ is incorporated into the Gaussian process predictive
variance, enabling the acquisition function to prioritize exploration in
regions exhibiting high local complexity.

\section{Sampling Allocation Theory for Inhomogeneous Bayesian Optimization}\label{sec:SATIBO}

This section develops a theoretical framework that explains how sampling effort should be distributed across regions with heterogeneous smoothness in BO. Building on the inhomogeneous BO framework introduced in Section~\ref{sec:IBO}, we analyze how local regularity properties of the objective function influence the decay of the $\mathbf{EI}$ acquisition function.  To make the theoretical analysis tractable, we assume throughout this section that the posterior mean and variance within each region $\X_i$ are constructed using only the local observation subset $\mathcal{D}_N^i$. 
In slight abuse of notation, we write the local means and variances as
\begin{equation*}
  m_{N_i}(x) = m_{N_i}(x ; \mathcal{D}^i_N) \quad \text{ and } \quad k_{N_i}(x) = k_{N_i}(x, x; \mathcal{D}_N^i) ,
\end{equation*}
where the subscript $N_i$ simultaneously denotes both the number of samples in $\X_i$ and that the local dataset $\mathcal{D}_N^i$ is being used.
Consequently, the regional Expected Improvement acquisition considered in this section is
\begin{equation} \label{eq:EIi2}
  \mathbf{EI}_i(x ; \mathcal{D}_N^i) = \big(m_{N_i}(x) - f(x_i^*)\big)\,\Phi\!\left(z_i(x)\right)
  + \hat{\sigma}_i \sqrt{\smash[b]{k_{N_i}(x)}} \,\phi\!\left(z_i(x)\right),
\end{equation}
if $k_{N_i}(x) > 0$, where $z_i(x) = [m_{N_i}(x) - f(x_i^*)] / (\hat{\sigma}_i \sqrt{\smash[b]{k_{N_i}(x)}})$ and $x_i^*$ is the best observed point in $\X_i$, in contrast to the best observed point $x^*$ in all of $\X$ used in Section~\ref{sec:IBO}.
Our goal is to characterize how differences in smoothness translate into non-uniform sampling and to derive principled allocation rules that balance exploration across regions.

\subsection{Expected Improvement Decay and Regional Allocation}\label{subsec:EIDRA}

We analyze the behavior of the $\mathbf{EI}_i(x; \mathcal{D}_N^i)$ under certain asymptotic assumptions on the local GP posterior. We show that the supremum of the acquisition function is governed by local decay rates. Consequently, ensuring that the $\mathbf{EI}$ is bounded by a common threshold across all subdomains can be achieved by allocating polynomially more evaluations to the less smooth regions.

For two positive sequences $u_N$ and $v_N$, we write $u_N \lesssim v_N$ if there exists a constant $C > 0$ such that $u_N \leq C v_N$ for all sufficiently large $N$. Similarly, we use $u_N \asymp v_N$ to indicate that $u_N$ and $v_N$ are of the same asymptotic order, meaning that both $u_N \lesssim v_N$ and $v_N \lesssim u_N$ hold.
\begin{ass} \label{ass1}
 The local GP posterior quantities within each subdomain $\X_i$ are assumed to satisfy the following asymptotic scaling properties with respect to the local sample size $N_i$ for each $i$:
    \begin{enumerate}
    \item Estimated local scale: $\hat{\sigma}_{i} \lesssim  N_i^{a_i}  $ with $a_i \in \R$. \label{ass11}
    \item Posterior standard deviation: \(\sup_{x\in \X_i}\sqrt{k_{N_i}(x)} \lesssim N_i^{-c}\) with $c -a_i > 0$. \label{ass12}
    \item Mean error: 
   $\sup_{x\in \X_i}\lvert m_{N_i}(x) - f(x)\rvert \lesssim N_i^{-b_i}$ with  $b_i > 0$. 
    \label{ass13}
    \item Regional regret: $\sup_{x\in \X_i}(f(x) - f(x_i^*)) \lesssim N_i^{-\rho_i}  $ with $ \rho_i > 0$. \label{ass14}
\end{enumerate}
\end{ass}

The four exponents describe three bottlenecks: $c-a_i$ measures how fast the \emph{scaled} posterior uncertainty contracts, $b_i$ the rate of convergence of the surrogate, and $\rho_i$ how fast the search closes the true optimality gap in $\X_i$. Larger $c - a_i$, $b_i$ and $\rho_i$ correspond to faster convergence and to smoother local behavior; they depend on both the kernel and the local regularity of $f$.

\begin{lem} \label{pro1}
  Suppose that Assumption~\ref{ass1} holds and let $\gamma_i \coloneqq \min(c - a_i,b_i,\rho_i)$.
  Then
       \begin{align}
           \sup_{x\in \X_i}\mathbf{EI}_i(x ; \mathcal{D}_{N}^i) \lesssim N_i^{-\gamma_i}.\label{eq:EIupper}
       \end{align}
\end{lem}
 \begin{proof}
Since $0\le\Phi\le1$ and $0\le\phi\le1$, equation~\eqref{eq:EIi2} gives for every $x\in\X_i$ with $k_{N_i}(x)>0$,
\[
  \mathbf{EI}_i(x ; \mathcal{D}_{N}^i)
  \le \max\big(0,\, m_{N_i}(x) - f(x_i^*)\big)
  + \hat{\sigma}_i \sqrt{\smash[b]{k_{N_i}(x)}} .
\]
The inequality is trivial when $k_{N_i}(x)=0$ because $\mathbf{EI}_i(x ; \mathcal{D}_{N}^i)=0$. Decomposing the mean gap
and applying Assumption~\ref{ass1}(\ref{ass13}) and~(\ref{ass14}) gives
\[
  m_{N_i}(x) - f(x_i^*)
  = \big(m_{N_i}(x)-f(x)\big) + \big(f(x)-f(x_i^*)\big)
  \lesssim N_i^{-b_i} + N_i^{-\rho_i}
  \lesssim N_i^{-\min(b_i,\rho_i)}.
\]
For the uncertainty term, Assumption~\ref{ass1}(\ref{ass11})
and~(\ref{ass12}) give
$\hat{\sigma}_i\sqrt{k_{N_i}(x)} \lesssim N_i^{a_i}N_i^{-c}
 = N_i^{-(c-a_i)}$.
Both bounds are uniform in $x\in\X_i$. Therefore
\[
  \sup_{x\in\X_i}\mathbf{EI}_i(x ; \mathcal{D}_{N}^i)
  \lesssim N_i^{-\min(b_i,\rho_i)} + N_i^{-(c-a_i)}
  \lesssim N_i^{-\min(c-a_i,\,b_i,\,\rho_i)}
  = N_i^{-\gamma_i}.
\]
\end{proof}

By Lemma~\ref{pro1} for each $i$ there exist a constant
$C_{i} > 0$ and a threshold $\underline{N}_i\in\N$ with 
\[
  \sup_{x\in\X_i}\mathbf{EI}_i(x ; \mathcal{D}_{N}^i)
  \le C_{i}\,N_i^{-\gamma_i}
  \quad \text{ for all } \quad N_i \ge \underline{N}_i .
\]
We fix one such pair for each region. For a target level $\delta>0$ set $N_i(\delta) = (C_i / \delta)^{1/\gamma_i}$.
$\gamma_i>0$, the map
$N \mapsto C_{i}N^{-\gamma_i}$ is strictly decreasing on $(0,\infty)$. Consequently,

\[
  \sup_{x\in\X_i}\mathbf{EI}_i(x ; \mathcal{D}_{N}^i) \le \delta
  \quad\text{ for every } \quad 
  N_i \geq \max\big(N_i(\delta),\, \underline{N}_i\big).
\]

\begin{thm}\label{th:thm0}
    Let the domain be partitioned into $\X_1$ and $\X_2$ and suppose that Assumption~\ref{ass1} holds with $\gamma_1 \le \gamma_2$, which holds in particular when $a_1 \ge a_2$, $b_1 \le b_2$ and $\rho_1 \le \rho_2$, so that $\X_1$ is the rougher subdomain. Set $r \coloneqq \gamma_2/\gamma_1$, $\kappa \coloneqq (C_1/C_2)^{1/\gamma_1}$, $\delta_0 \coloneqq \min_{i\in \{1,2\}} C_i\underline{N}_i^{-\gamma_i}$ and $\delta \in (0,\delta_0]$.   
    If
    \begin{equation}\label{eq:alloc_cond}
    N_2 \geq N_2(\delta)
    \qquad\text{and}\qquad
    N_1 \geq \kappa N_2^{r},
    \end{equation}
    then $\sup_{x\in\X_i}\mathbf{EI}_i(x;\mathcal{D}_N^i)\le\delta$ for $i \in \{1, 2\}$. 

\end{thm}
\begin{proof}
Since $\delta\le\delta_0$, we have $ N_i(\delta)\ge \underline{N}_i$ for $i \in \{1,2\}$. Therefore, whenever $N_i \ge N_i(\delta)$, the bound of Lemma~\ref{pro1} is applicable and gives $\sup_{x\in\X_i}\mathbf{EI}_i(x;\mathcal{D}_N^i)\le C_i N_i^{-\gamma_i}$.
Suppose \eqref{eq:alloc_cond} holds. The first inequality gives
$C_2N_2^{-\gamma_2} \le \delta$ directly. For the first region, $N_2 \ge N_2(\delta)$ and $r>0$ give $N_2^{r} \ge N_2(\delta)^{r}$, so the second
inequality yields
\[
  N_1 \ge \kappa N_2^{ r} \ge \kappa N_2(\delta)^{ r}
  = \Big(\frac{C_1}{C_2}\Big)^{1/\gamma_1}
        \Big(\frac{C_2}{\delta}\Big)^{\gamma_2/(\gamma_2\gamma_1)}
  = \Big(\frac{C_1}{\delta}\Big)^{1/\gamma_1}
  = N_1(\delta) .
\]
Hence $C_1N_1^{-\gamma_1}\le\delta$. Both suprema are therefore bounded by
$\delta$.
\end{proof}

Theorem~\ref{th:thm0} quantifies the imbalance due to heterogeneous smoothness: ensuring that the maximum Expected Improvement falls below a shared target level $\delta$ across both regions is achieved by an allocation that places polynomially more evaluations in the rougher subdomain $\X_1$ than in $\X_2$, with an exponent $r=\gamma_2/\gamma_1$ determined by the ratio of the local decay rates. Equivalently, if the local sample sizes grow at comparable rates, the bound of Lemma~\ref{pro1} on the maximum Expected Improvement decays faster in the smoother region, so the certification of a common level is limited by the rougher subdomain.

\begin{cor}\label{th:thm1}
Let $\X$ be partitioned into $n\ge2$ subdomains $\{\X_l\}_{l=1}^{n}$ and suppose that Assumption~\ref{ass1} holds.
Let the subdomains be indexed so that $\gamma_1\le\dots\le\gamma_n$, which
holds in particular when $a_1\ge\dots\ge a_n$, $b_1\le\dots\le b_n$ and
$\rho_1\le\dots\le\rho_n$. For $i\ne j$ set
$r_{ij}\coloneqq\gamma_j/\gamma_i$ and
$\kappa_{ij}\coloneqq(C_i/C_j)^{1/\gamma_i}$, and let
$\delta_0\coloneqq\min_l C_l\underline N_l^{-\gamma_l}$ and
$\delta\in(0,\delta_0]$. If
\[
  N_n\ge N_n(\delta)
  \qquad\text{and}\qquad
  N_l\ge \kappa_{ln}\,N_n^{\,r_{ln}}
  \quad\text{ for } \quad l \in \{1,\dots,n-1\},
\]
then $\sup_{x\in\X_l}\mathbf{EI}_l(x ; \mathcal{D}_{N}^l)\le\delta$ for every $l$. 
\end{cor}

\begin{proof}
Lemma~\ref{pro1} involves only the subdomain index, so the argument of
Theorem~\ref{th:thm0} applies verbatim to each pair $(\X_l,\X_n)$; no interaction between three or more subdomains arises. Since the
hypothesis constrains every subdomain through its pairing with $\X_n$, the
conclusion holds for all $l$. 
\end{proof}

The scaling assumptions can be naturally extended to allow the posterior uncertainty and mean error to depend on the specific local kernel $k_{i}$. Unlike the single-kernel setting where the uncertainty decay rate $c$ is a global constant, incorporating subdomain-specific rates $c_{i}$ reflects that different kernel families converge at different rates.  
By adapting Assumption~\ref{ass1} to incorporate these local rates $c_{i}$, Corollary~\ref{th:thm1} holds verbatim with $\gamma_i = \min(c_i - a_i, b_i, \rho_i)$.
The ordering $\gamma_1 \leq \cdots \leq \gamma_n$ holds in particular when $c_i - a_i$, $b_i$ and $\rho_i$ are all non-decreasing in $i$.
The conclusion is the same: the sampling effort sufficient to certify a common acquisition ceiling grows polynomially faster in regions of lower smoothness, with exponent $r_{ij}=\gamma_j/\gamma_i \ge 1$ for $i < j$, strictly greater than one when $\gamma_i < \gamma_j$. This principle forms the foundation of IBO.

\subsection{Sampling Allocation under Sobolev Smoothness} \label{subsec:sobo}

We now ground the scaling results of Section~\ref{subsec:EIDRA} in a concrete functional-analytic setting by assuming that the objective function exhibits piecewise Sobolev regularity and that the surrogate model is based on kernels with a known Reproducing Kernel Hilbert Space (RKHS).
To establish the necessary function spaces for our analysis, we introduce Sobolev Hilbert spaces via the Fourier transform. 
The Fourier transform is $\hat{f}(\xi) = \int_{\mathbb{R}^d} f(x) e^{-2\pi i x \cdot \xi} dx$.

\begin{defn}
    Let $\alpha \ge 0$ be a real number. The (fractional) Sobolev space $\Hi^\alpha(\R^d)$ is defined as the space of functions $f \in L^2(\R^d)$ whose Fourier transform $\hat{f}$ satisfies
    $$\|f\|_{\Hi^\alpha(\R^d)} = \left( \int_{\R^d} (1 + \|\xi\|^2)^\alpha |\hat{f}(\xi)|^2 \, d\xi \right)^{1/2} < \infty.$$
    Equipped with the inner product 
    $$\langle f, g \rangle_{\Hi^\alpha(\R^d)} = \int_{\R^d} (1 + \|\xi\|^2)^\alpha \hat{f}(\xi) \overline{\hat{g}(\xi)} \, d\xi,$$
    the space $\Hi^\alpha(\R^d)$ forms a Hilbert space.
    For a domain $\X \subset \R^d$, the space $\Hi^{\alpha}(\X)$ is defined as the restriction of $\Hi^{\alpha}(\R^d)$ to $\X$ equipped with the quotient norm $$\|f\|_{\Hi^{\alpha}(\X)} = \inf \{ \|g\|_{\Hi^{\alpha}(\R^d)} : g|_{\X} = f \}.$$
\end{defn}
When $\alpha$ is a non-negative integer, this Fourier-analytic definition is norm-equivalent to the classical Sobolev space $W^{\alpha,2}(\R^d) = \{ f \in L^2(\R^d) : D^s f \in L^2(\R^d) \text{ for all } |s| \leq \alpha \}$ equipped with the norm $\|f\|_{W^{\alpha,2}(\R^d)} = ( \sum_{|s| \leq \alpha} \|D^s f\|_{L^2(\R^d)}^2 )^{1/2}$, where $D^s f$ denotes the weak derivative~\cite{di2012hitchhikers}.

\begin{defn}
    Let $\X$ be a non-empty set and $\Hi$ be a Hilbert space of real-valued functions defined on $\X$, equipped with an inner product $\langle \cdot, \cdot \rangle_{\Hi}$. The space $\Hi$ is an RKHS if the evaluation functional $f \mapsto f(x)$ is continuous for all $x \in \X$. 
\end{defn}

By the Riesz Representation Theorem, this continuity implies that for every $x \in  \X$ there exists a unique element $k_x \in \Hi$ such that for every $f \in \Hi$, the \emph{reproducing property} holds: $f(x) = \langle f, k_x \rangle_{\Hi}$. By defining the function $k(x, x') = k_x(x') = \langle k_x, k_{x'} \rangle_{\Hi}$, it follows that $k$ is a symmetric, positive definite kernel termed the \emph{reproducing kernel}. This kernel $k$ completely characterizes the space, which is commonly denoted by $\Hi(k)$. 
Conversely, the Moore--Aronszajn theorem~\cite[Section~2, p.~344]{aronszajn1950theory} establishes that every symmetric, positive definite kernel on $ \X \times  \X$ uniquely defines an RKHS for which it is the reproducing kernel. 
A key result in the theory of RKHS is the equivalence between Mat\'{e}rn spaces and Sobolev spaces. When the kernel $k$ is a Mat\'{e}rn kernel \eqref{eq:matern} with smoothness $\nu$, the resulting RKHS, $\Hi(k)$, is norm-equivalent to the Sobolev space $\Hi^{\alpha}(\R^d)$ for $\alpha = \nu + d/2$.

\begin{ass}\label{ass2}
 Assume that the following hold:
 \begin{enumerate}
    \item  RKHS: The kernel $k$ induces an RKHS, denoted $\Hi(k)$, that is norm-equivalent to the Sobolev space $\Hi^{\alpha}(\X)$ with $\alpha > d/2$, which ensures that $\Hi^{\alpha}(\X)$ embeds in $C(\X)$. \label{ass21}
    \item Domain regularity: Each subdomain $\X_i \subset \R^d$ is compact, connected, has a non-empty interior, and satisfies an interior cone condition with a Lipschitz boundary (e.g., \cite[Theorem 5.24]{adams2003sobolev}). These conditions ensure the existence of a bounded Sobolev extension operator from $\Hi^\alpha(\X_i)$ to $\Hi^\alpha(\R^d)$ for every $\alpha>0$. \label{ass22}
    \item  Local regularity: For each $i \in \{1, \dots, n\}$ there are $\beta_i \in (d/2,\alpha]$
    and $B_i < \infty$ such that $\lVert f|_{\X_i} \rVert_{\Hi^{\beta_i}(\X_i)} \le B_i$. \label{ass23} 
    \item Quasi-uniformity: For each $i \in \{1, \dots, n\}$, the design  $X_i$ is quasi-uniform in $\X_i$. That is, defining the fill distance and the separation radius respectively as
    \begin{equation*}
      h_{X_i,\X_i} = \sup_{x\in\X_i}\min_{x_j\in X_i}\lVert x-x_j\rVert
    \qquad \text{and} \qquad
    q_{X_i} = \tfrac12\min_{x_j\neq x_l\in X_i}\lVert x_j-x_l\rVert,
    \end{equation*}
    the local mesh ratios
    $h_{X_i,\X_i}/q_{X_i}$ remain bounded as $N_i\to\infty$ and
    $h_{X_i,\X_i}\asymp N_i^{-1/d}$.
    \label{ass24}
\end{enumerate}
\end{ass}

These assumptions allow the model to misspecify the smoothness of $f$.
Specifically, whenever $\beta_i < \alpha$, the objective is rougher than the kernel assumes and need not lie in the global RKHS $\Hi(k)$. The differences between the local regularities $\beta_i$ drive the allocation law of Theorem~\ref{cor:cor42}. 
The local norm bounds $B_i$ enter only through the hidden constants and do not affect the rates of convergence.
Within the context of this assumption, Karvonen et al.~\cite{karvonen2020maximum} establish asymptotic results for the ML estimates of $\sigma$, the posterior variance and the absolute error of the GP estimate.\footnote{In our setting the assumption in~\cite{karvonen2020maximum} that $\lfloor \beta_i \rfloor > d/2$ can be relaxed to $\beta_i > d/2$; see the discussion in Section~4.3 of~\cite{karvonen2020maximum}.}
These asymptotics play crucial roles in the proofs below.

\begin{lem}[Regret upper bound]\label{pro:regret_upper}
Suppose that Assumption~\ref{ass2} holds and let $x_i^\dagger \in \argmax_{x\in\X_i} f(x)$. 
Set $p_i = 2$ if $x_i^\dagger \in \operatorname{int}\X_i$ and $p_i = 1$ otherwise.  Define $s_i \coloneqq \min(\beta_i - d/2, p_i)$ and $\rho_i \coloneqq s_i/d$. 
If $\beta_i - d/2 \notin \{1, p_i\}$, then
\begin{equation}\label{eq:regret}
  \sup_{x\in\X_i}\big(f(x)-f(x_i^{*})\big) \lesssim N_i^{-\rho_i}.
\end{equation}
In particular, $\rho_i = \beta_i/d - 1/2$ whenever $\beta_i < d/2+p_i$.
\end{lem}

\begin{proof}
Since $\beta_i>d/2$, Sobolev embedding makes $f$ continuous on $\X_i$, which is
compact by Assumption~\ref{ass2}(\ref{ass22}), so the maximizer $x_i^\dagger$
exists. By Assumption~\ref{ass2}(\ref{ass24}) the design $X_i$ has fill
distance $h_{X_i,\X_i} \asymp N_i^{-1/d}$, so there are points $x_j \in X_i$ such that
\begin{equation}\label{eq:close}
  \lVert x_i^\dagger - x_j \rVert \le h_{X_i,\X_i} \asymp N_i^{-1/d}.
\end{equation}
Because $x_i^{*}$ is the current optimal point in $\X_i$, $f(x_j) \le f(x_i^{*})$ and therefore
\begin{equation}\label{eq:reduce}
  \sup_{x\in\X_i}\big(f(x)-f(x_i^{*})\big)
  = f(x_i^\dagger)-f(x_i^{*})
  \le f(x_i^\dagger)-f(x_j).
\end{equation}
Next we show that for every $y\in\X_i$ sufficiently close to $x_i^\dagger$,
\begin{equation}\label{eq:increment}
  f(x_i^\dagger)-f(y) \lesssim \lVert y-x_i^\dagger\rVert^{\,s_i}.
\end{equation}

If $\beta_i-d/2<1$, then $s_i=\beta_i-d/2<1$ regardless of $p_i$. The Sobolev embedding theorem gives $f|_{\X_i}\in C^{0,s_i}(\X_i)$ with H\"older seminorm bounded by a multiple of $B_i$, which is \eqref{eq:increment}.
If $\beta_i-d/2>1$ and $p_i=1$, then $s_i=1$. The Sobolev embedding theorem gives $f|_{\X_i}\in C^{0,1}(\X_i)$ with Lipschitz seminorm bounded by a
multiple of $B_i$, which is \eqref{eq:increment} with $s_i = 1$. The maximizer plays no part in this branch, so the bound holds regardless of  $x_i^\dagger$ being in the interior of $\X_i$ or not.
If $\beta_i-d/2>1$ and $p_i=2$, then each component of $\nabla f$ lies in
$\Hi^{\beta_i-1}(\X_i)$ with $\beta_i-1>d/2$, so the Sobolev embedding theorem
gives $\nabla f \in C^{0,\theta}(\X_i)$ with
$\theta = \min(\beta_i-d/2-1,\,1) = s_i-1$.
In particular, $f|_{\X_i} \in C^{1}(\X_i)$ and $\nabla f(x_i^\dagger)=0$, because
$x_i^\dagger$ is an interior maximizer. Writing $v=y-x_i^\dagger$ and
integrating along the segment from $x_i^\dagger$ to $y$, which lies in $\X_i$
once $\lVert v\rVert$ is small because $x_i^\dagger$ is in the interior of $\X_i$,
\[
  \big\lvert f(y)-f(x_i^\dagger)\big\rvert
  = \Big\lvert \int_0^1 \big(\nabla f(x_i^\dagger+tv)-\nabla f(x_i^\dagger)\big)
    \cdot v \,dt \Big\rvert
  \lesssim \lVert v\rVert^{\,1+\theta},
\]
and $1+\theta=s_i$, which is again \eqref{eq:increment}.

Taking $y=x_j$ in \eqref{eq:increment}, which is admissible for $N_i$ large by
\eqref{eq:close}, and combining with \eqref{eq:reduce} and \eqref{eq:close} yields
\[
  \sup_{x\in\X_i}\big(f(x)-f(x_i^{*})\big)
  \le f(x_i^\dagger)-f(x_j)
  \lesssim \lVert x_i^\dagger-x_j\rVert^{\,s_i}
  \lesssim N_i^{-s_i/d}
  = N_i^{-\rho_i}.
\]
\end{proof}

\begin{lem} \label{cor:cor41}
Suppose that Assumption~\ref{ass2} holds and $\beta_i-d/2 \notin \{1, p_i\}$, where $p_i$ and $\rho_i$ are defined as in Lemma~\ref{pro:regret_upper}. Then $\gamma_i = \rho_i$ and
\[
\sup_{x\in \X_i}\mathbf{EI}_i(x ; \mathcal{D}_{N}^i) \lesssim N_i^{-\rho_i }.
\]
\end{lem}

\begin{proof}
We use the rates cited and derived in~\cite{karvonen2020maximum} to show that Assumption~\ref{ass2} implies Assumption~\ref{ass1}. 
Theorems~4.2 and~4.4 and Proposition~4.5 in \cite{karvonen2020maximum} yield
\begin{equation*}
    \sup_{x\in \X_i} \sqrt{\smash[b]{k_{N_i}(x)}} \lesssim N_i^{-\frac{\alpha}{d}+\frac{1}{2}}, \quad \sup_{x\in \X_i} \lvert m_{N_i}(x)-f(x)\rvert \lesssim N_i^{-\frac{\beta_i}{d}+\frac{1}{2}} \quad \text{ and } \quad \hat{\sigma}_{i} \lesssim N_i^{\frac{\alpha-\beta_i}{d}-\frac{1}{2}} .
\end{equation*}
In the notation of Assumption~\ref{ass1}, 
\begin{equation*}
  c = \alpha/d - 1/2, \quad a_i = (\alpha - \beta_i) / d - 1/2, \quad \text{ and } \quad b_i = \beta_i / d - 1 / 2.
\end{equation*}

Therefore $c - a_i = \beta_i / d > b_i$.
By Lemma~\ref{pro:regret_upper}, the regional regret decays at the rate $\rho_i = \min(\beta_i/d - 1/2,p_i/d) = \min(b_i, p_i/d) \leq b_i$.
Therefore
\begin{equation*}
  \sup_{x\in \X_i}\mathbf{EI}_i(x ; \mathcal{D}_{N}^i) \lesssim N_i^{-\gamma_i }
\end{equation*}
with $\gamma_i = \min(c - a_i, b_i, \rho_i) = \rho_i$ in Lemma~\ref{pro1}. 
\end{proof}

\begin{thm}\label{cor:cor42}
Suppose that Assumption~\ref{ass2} holds.
Let $\X$ be partitioned into $n\ge2$ subdomains $\{\X_k\}_{k=1}^{n}$ indexed such that
\(
  \beta_1 \le \beta_2 \le \dots \le \beta_n
\) and $\rho_1 \le \rho_2 \le \dots \le \rho_n$.
Suppose $\beta_k - d/2 \notin \{1,p_k\}$ for every $k$, with $p_k$ as in Lemma~\ref{pro:regret_upper}. Then $\gamma_k = \rho_k$ for every $k$, and for any pair
$i < j$ the allocation exponent of Corollary~\ref{th:thm1} is
\[
  r_{ij} = \frac{\rho_j}{\rho_i}
  = \frac{\min\big(2\beta_j - d, 2p_j\big)}
             {\min\big(2\beta_i - d, 2p_i\big)}.
\]
\end{thm}

\begin{proof}
By Lemma~\ref{cor:cor41}, the regime is regret-dominant, so
$\gamma_k=\min(c-a_k,b_k,\rho_k)=\rho_k$, and
Lemma~\ref{pro:regret_upper} gives $\rho_k= \min (\beta_k/d-1/2,p_k/d)$. Hence
\[
  r_{ij}=\frac{\gamma_j}{\gamma_i}=\frac{\min (\beta_j/d-1/2,p_j/d)}{\min (\beta_i/d-1/2,p_i/d)}
        =\frac{\min \big(2\beta_j - d,2p_j\big)}{\min \big(2\beta_i - d,2p_i\big)}.
\]
\end{proof}

If $p_i = p_j$, the exponent satisfies $r_{ij} > 1$ only when the rougher subdomain is unsaturated, meaning $\beta_i < d/2+p_i$.
If $\beta_i \ge d/2+p_i$ and $\beta_j \ge d/2+p_j$ then $r_{ij}=p_j/p_i = 1$, so $r_{ij}>1$ requires $\beta_i<\beta_j$ together with $\beta_i<d/2+p_i$.

\section{Numerical Illustration}\label{sec:NIL}

Throughout this section, all benchmark evaluations are framed as minimization problems, which is equivalent to maximizing $-f$. The Python code to reproduce the experiments is available online\footnote{\url{https://github.com/RutherDondj/IBO-Framework}}.
 We evaluate the proposed framework from Section~\ref{sec:IBO} across a range of benchmark settings designed to assess its performance under different levels of complexity.
The evaluation is structured in three parts. We begin with standard homogeneous benchmark functions to establish a baseline comparison against Classic BO. We then consider synthetic inhomogeneous functions generated from non-stationary Gaussian processes to examine performance in spatially varying landscapes. Finally, we assess the method on a real-world 8-dimensional optimization task involving LASSO hyperparameter tuning.

To ensure a fair comparison, both Classic BO and all IBO variants employ the same surrogate kernel across all experiments: a Matérn kernel with smoothness parameter $\nu = 5/2$. Before any optimization run, the objective is evaluated at an independent Sobol design and the length-scale $\ell$ is estimated by maximum marginal likelihood~\cite[Chapter~5]{williams2006gaussian} on those evaluations, then held fixed throughout the run.

\subsection{Evaluation Metrics}

To ensure a robust comparison across different acquisition functions $\text{acq} \in \{\text{EI, PI, LCB, TS}\}$ and methods $m$, we define the following metrics used throughout all experiments.
Because the quality of the optimization heavily depends on the initial design, we conduct $n_{\text{seeds}} = 30$ independent runs to ensure the results are not biased by a single lucky or unlucky starting configuration. In each run, the randomness stems from the uniform selection of starting points. By averaging over these seeds, we isolate the true performance of the method from the variance introduced by the initial sampling.

\textbf{Mean Convergence Curve:} Let $f_{s}^{(m, \text{acq})}(t)$ denote the best objective value found up to iteration $t$ for a specific seed $s$. The mean convergence curve is defined as
\begin{align}\label{eq:mu}
    \mu(m, \text{acq})(t) = \frac{1}{n_{\text{seeds}}} \sum_{s=1}^{n_{\text{seeds}}} f_s^{(m, \text{acq})}(t).
\end{align}

\textbf{Aggregated Performance Metric:} To pool benchmarks of different scales, we
measure progress relative to the optimum. Every benchmark is shifted so that
its minimum is $f_{\min}=1$, and we define the normalized regret
\begin{equation}\label{eq:nr}
  \mathrm{NR}^{(m,\mathrm{acq})}_{p,s}(t)
  =
  \frac{f^{(m,\mathrm{acq})}_{p,s}(t) - f_{\min}}
       {f^{(m,\mathrm{acq})}_{p,s}(0) - f_{\min}} ,
\end{equation}
which starts at $1$ and decreases towards $0$, so a value of $0.05$ means that
$5\%$ of the initial error remains. All methods share the same initial design
within a seed, so the denominator is common to them and the quality of that design affects them all equally.

Averaging over seeds and then over the $n_f$ benchmark functions gives
the global convergence curve regret
\begin{equation}\label{eq:global_nr}
  \overline{\mathrm{NR}}^{(m,\mathrm{acq})}(t)
  =
  \frac{1}{n_f}\sum_{p=1}^{n_f}
  \frac{1}{n_{\mathrm{seeds}}}\sum_{s=1}^{n_{\mathrm{seeds}}}
  \mathrm{NR}^{(m,\mathrm{acq})}_{p,s}(t) .
\end{equation}

\subsection{Performance on Homogeneous Benchmark Functions}

This subsection evaluates the proposed IBO against Classic BO across a suite of standard homogeneous benchmark functions. We compare four configurations: Classic BO, $\text{IBO}_1$ (no partition with global $\sigma$ estimation), $\text{IBO}_4$ (four cells, halving along the first two coordinates), and $\text{IBO}_P$ (adaptive tree-based partitioning). Based on the acquisition framework defined in Section~\ref{subsec:AC}, we evaluate the convergence behavior of these methods across four acquisition functions: Expected Improvement (EI), Probability of Improvement (PI), Lower Confidence Bound (LCB), and Thompson Sampling (TS).

\subsubsection{Convergence Results: Homogeneous Functions}

\begin{figure}[t]
  \centering
  
  \subfloat[Kim1 function\label{fig:f_Eggholder}]{
    \includegraphics[width=0.48\textwidth]{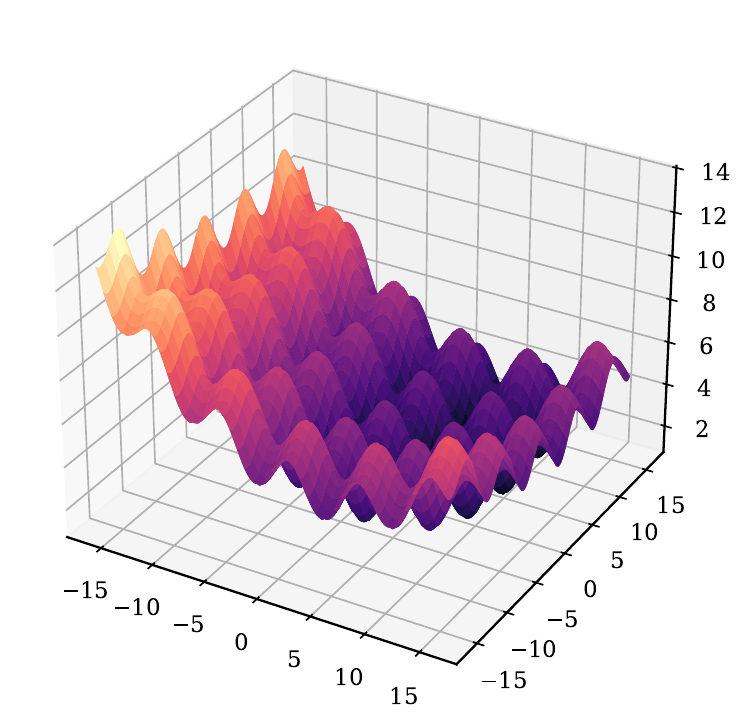}
  }
  \hfill
  \subfloat[Synthetic inhomogeneous function\label{fig:f_nstat}]{
    \includegraphics[width=0.48\textwidth]{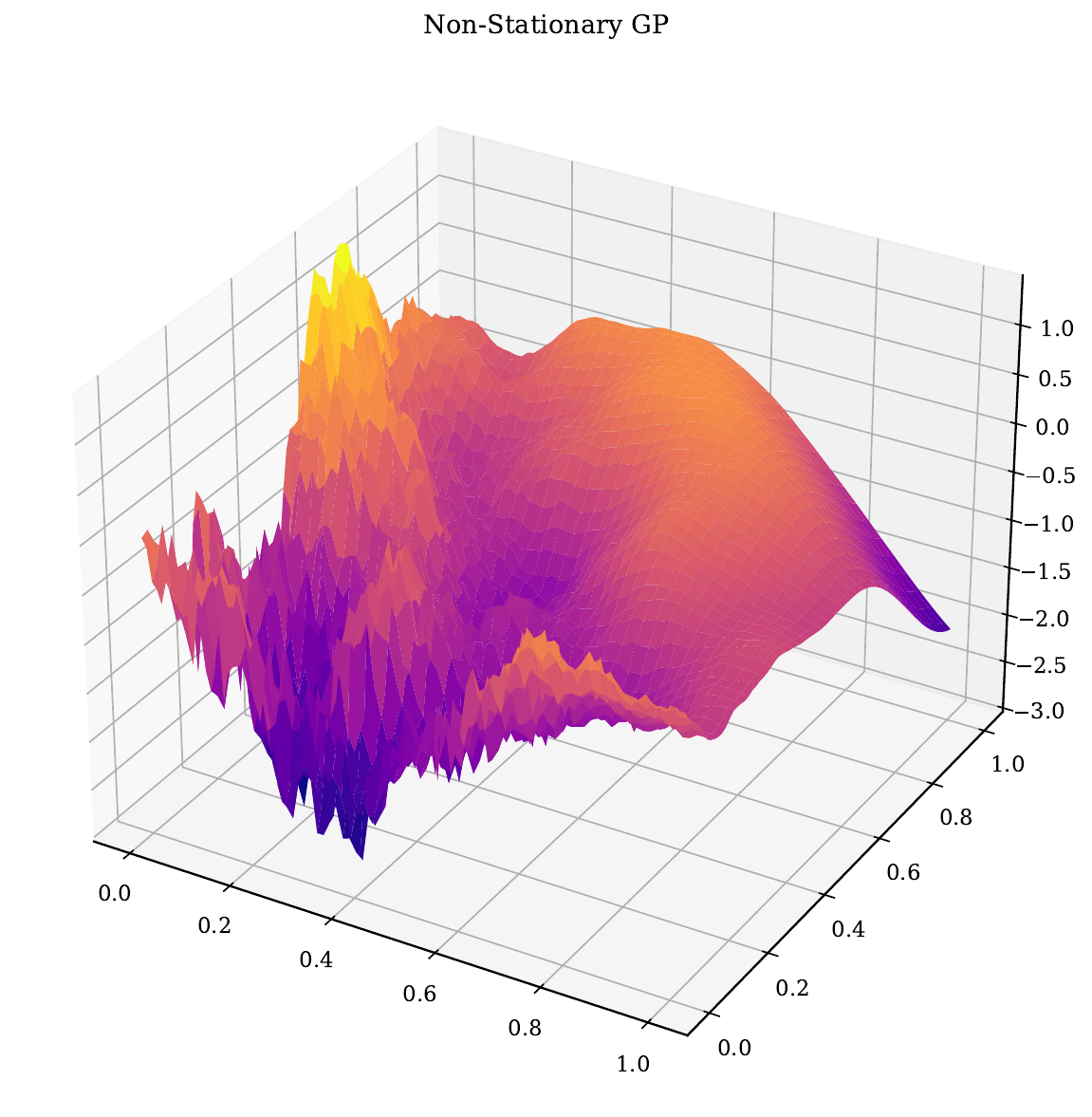}
  }

  \caption{(a) The Kim1 function exhibiting complex multi-modal topology; 
  (b) realization from an inhomogeneous Gaussian process.}
  
  \label{fig:f_mix}
\end{figure}

The Kim1 function~\cite[Eq.~27]{kim2020benchmark} of Figure~\ref{fig:f_Eggholder} superposes a sinusoidal ripple on an asymmetric quadratic bowl, and is therefore non-convex with numerous local minima.

\begin{figure}[t]
\centering
\includegraphics[width=\textwidth]{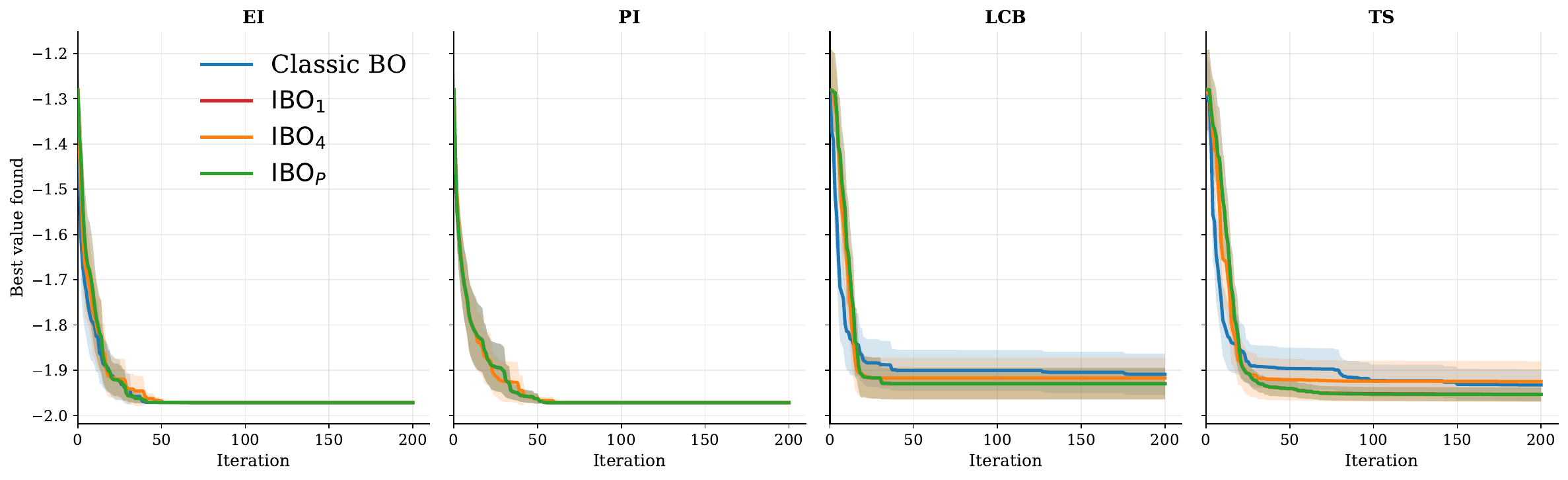}
\caption{Mean convergence curves for the Kim1 function across acquisition functions and methods. Shaded areas represent the 95\% confidence interval across seeds.}
\label{fig:eggholder_conv}
\end{figure}

In terms of exploration efficiency, as illustrated in Figure~\ref{fig:eggholder_conv}, which reports the mean convergence curve~\eqref{eq:mu}, the IBO methods outperform Classic BO in the LCB and TS variants. While Classic BO plateaus early (approx.\ iteration 25), $\text{IBO}_P$ and $\text{IBO}_4$ continue to improve before settling at a lower value, suggesting that localizing the scale $\sigma$ prevents the global GP from over-smoothing high-variance regions.
$\text{IBO}_P$ (green) consistently matches or exceeds the performance of fixed partitioning, demonstrating that adaptive splitting effectively identifies regions requiring higher model flexibility without manual grid tuning.

\subsubsection{Aggregate Results}

\begin{figure}[t]
\centering
\includegraphics[width=\textwidth]{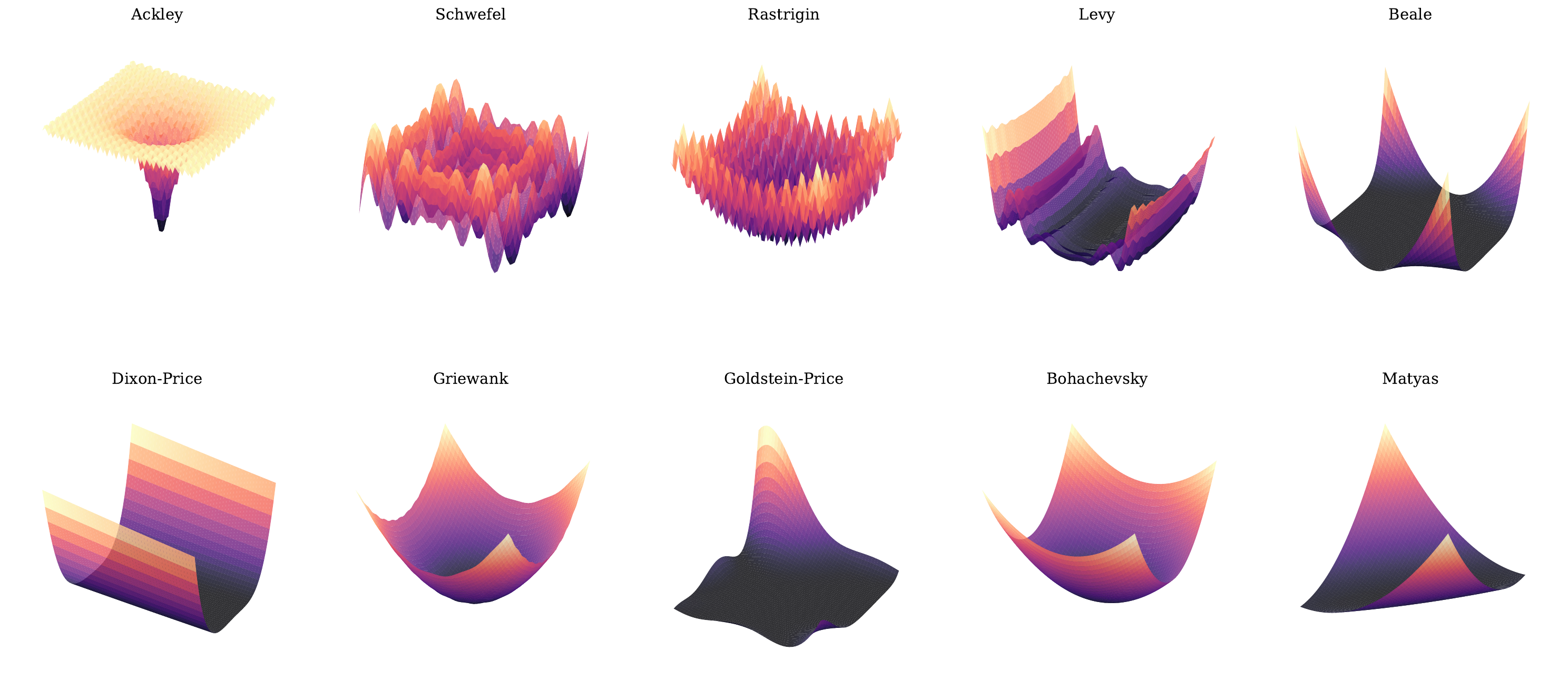}
\caption{$10$ benchmark homogeneous functions.}
\label{fig:10_benc}
\end{figure}

To evaluate the generalizability of the proposed framework, we aggregate results across $10$ standard homogeneous benchmark functions, including the Ackley, Schwefel, and Rastrigin surfaces (see Figure~\ref{fig:10_benc}).

\begin{figure}[t]
\centering
\includegraphics[width=\textwidth]{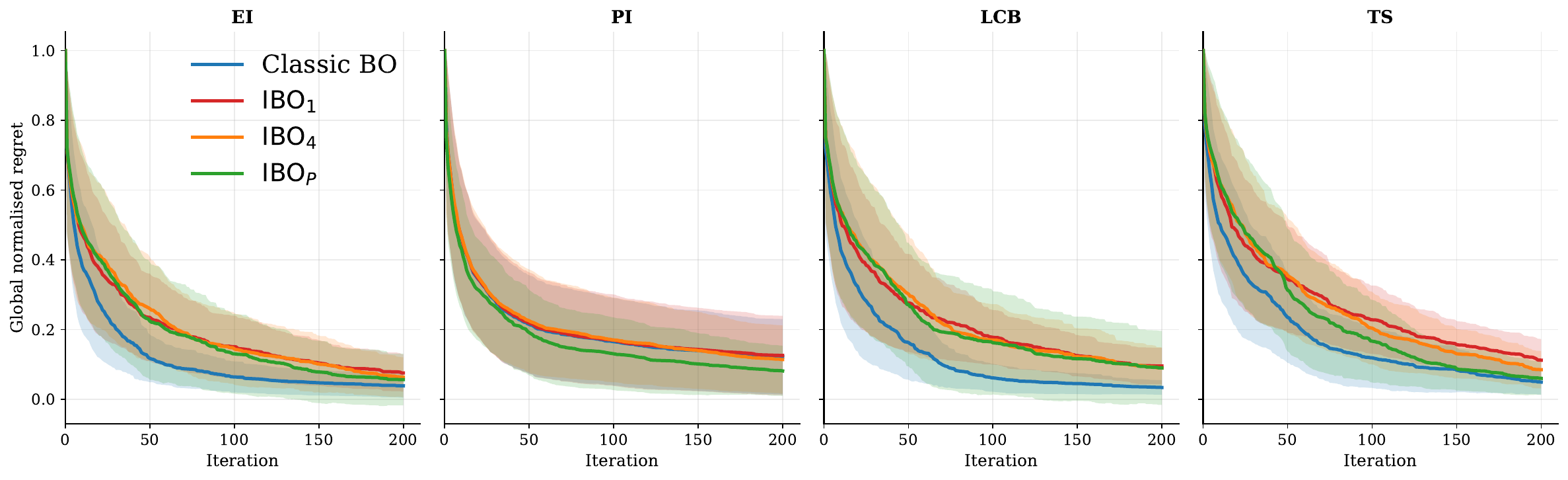}
\caption{Global normalized regret on the $10$ homogeneous functions averaged over the functions. Shaded regions are $95\%$ confidence intervals across functions.}
\label{fig:bench_ratios}
\end{figure}

Figure~\ref{fig:bench_ratios} reports the global normalized regret~\eqref{eq:global_nr} on the $10$ homogeneous functions. Classic BO converges faster under EI, LCB and TS, while $\text{IBO}_P$ outperforms it under PI; the confidence bands overlap throughout. 
The final values are closer than the curves suggest: Table~\ref{tab:main} gives a mean difference against Classic BO that is significant only under PI, in favor of $\text{IBO}_P$, and under LCB, against it, and Table~\ref{tab:ablation} puts $\text{IBO}_P$ ahead on $23$ of the $40$ function--acquisition pairs with a median ratio of $1.12$.
The criterion partitions these objectives in $68\%$ of runs, with $4.36$ leaves on average, because stationarity of the objective does not imply a spatially constant scale: Levy, Beale, Dixon-Price
and Goldstein-Price span several orders of magnitude and are partitioned in every run, whereas Bohachevsky and Matyas are never partitioned.

\subsection{Performance on Inhomogeneous Functions}\label{sec:nonstat}

While the previous evaluation established a baseline on homogeneous functions, the primary utility of the IBO framework lies in its ability to adapt to inhomogeneous ones.
To test this adaptability, we generate objective functions from a non-stationary GP, $f \sim \mathcal{GP}(0, K_{NS})$, constructed as
\[
  f(x)=\sum_{i=1}^{n}\tilde w_i(x)\,g_i(x),
  \qquad g_i\sim\mathcal{GP}(0,k_i)\ \text{independent}.
\]
By independence, $K_{NS}(x, x') = \sum_{i=1}^{n} \tilde{w}_i(x) \tilde{w}_i(x') k_i(x, x')$, which is a valid covariance.

The kernels $k_i$ are stationary with a \emph{common} variance, $k_i(x,x)=\sigma^2$ for all $i$.
The base weights
$w_i(x)=\exp\{-\|x-\boldsymbol{\mu}_i\|^2/(2\ell_w^2)\}$
are Gaussian bumps centered at anchor points $\boldsymbol{\mu}_i$.
The weights $\tilde{w}_i$ are obtained via $\ell^2$-normalization:
\[
  \tilde w_i(x)=\frac{w_i(x)}{\sqrt{\sum_{j=1}^{n}w_j(x)^2}} .
\]
Consequently, the pointwise variance is constant: $K_{NS}(x,x)=\sigma^2\sum_i\tilde w_i(x)^2=\sigma^2$ for all $x\in \X$. The spatial zones therefore differ only in their correlation structure, not their marginal amplitude. Since the strictly positive Gaussian weights decay to near-zero away from their anchor points, the zones differ in their \emph{effective} roughness rather than their Sobolev regularity.

To generate realizations from the non-stationary GP $f$, discrete samples $\mathbf{y}$ evaluated over a spatial grid are obtained via the standard Cholesky decomposition of the synthesized covariance matrix $\mathbf{K} = \mathbf{L}\mathbf{L}^\top$, such that $\mathbf{y} = \mathbf{L}\mathbf{z}$ with $\mathbf{z} \sim \mathcal{N}(\mathbf{0}, \mathbf{I})$~\cite{williams2006gaussian}.

\subsubsection{Convergence Results: Inhomogeneous Function}

We first analyze a specific non-stationary realization generated with three fixed anchor points located at $[0.2, 0.2]$ (rough zone), $[0.5, 0.5]$ (transitional zone), and $[0.8, 0.8]$ (smooth zone), using a weight length-scale $\ell^2_w = 0.05$ to control the transition between them. These zones are governed by basis kernels $k_i$ representing three distinct levels of smoothness:
\begin{itemize}
    \item Rough ($\nu=1/2$, $\ell = 0.2$): Mat\'{e}rn $k_{1}(r) = \sigma^2 \exp(-r/\ell)$
    \item Transitional ($\nu=3/2$, $\ell = 0.4$): Mat\'{e}rn $k_{2}(r) = \sigma^2 (1 + \frac{\sqrt{3}r}{\ell}) \exp(-\frac{\sqrt{3}r}{\ell})$
    \item Smooth ($\nu = 5/2$, $\ell = 0.3$): Mat\'{e}rn $k_{3}(r) = \sigma^2 \big(1 + \frac{\sqrt{5}r}{\ell}
    + \frac{5r^2}{3\ell^2}\big)\exp\big(-\frac{\sqrt{5}r}{\ell}\big)$
\end{itemize}
The resulting surface, characterized by spatially varying effective roughness, is depicted in Figure~\ref{fig:f_nstat}.

\begin{figure}[t]
\centering
\includegraphics[width=\textwidth]{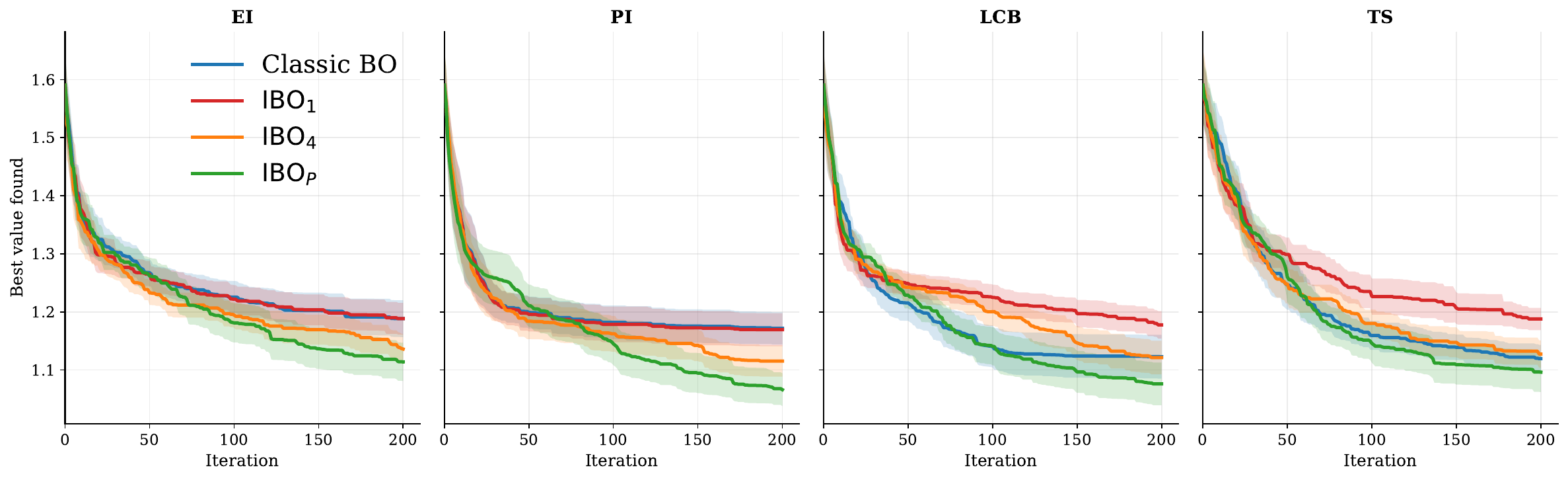}
\caption{Mean convergence curves on the inhomogeneous function of Figure~\ref{fig:f_nstat}. Shaded areas represent the 95\% confidence interval across seeds.}
\label{fig:nstat_conv}
\vspace{-1pt} 
\end{figure}

The convergence curves for this surface (Figure~\ref{fig:nstat_conv}), which reports the mean convergence curve~\eqref{eq:mu}, reveal a more pronounced divergence in performance compared to the homogeneous case.
A partitioning advantage is apparent across all four acquisition functions. Most notably in EI and PI, the partitioned methods ($\text{IBO}_{4}$ and $\text{IBO}_{P}$) outperform Classic BO and $\text{IBO}_{1}$. 
Moreover, $\text{IBO}_{P}$ (green) attains the smallest final value under every acquisition function. By recursively splitting the domain when the non-stationarity ratio $\eta$ exceeds the threshold $\tau$, the tree structure assigns more model flexibility to the rough regions while maintaining broader coverage in smooth regions. This results in superior regret reduction compared to $\text{IBO}_{4}$, which uses a fixed partition that may not align with the local landscape features.

\subsubsection{Aggregate Results}

\begin{figure}[t]
\centering
\includegraphics[width=\textwidth]{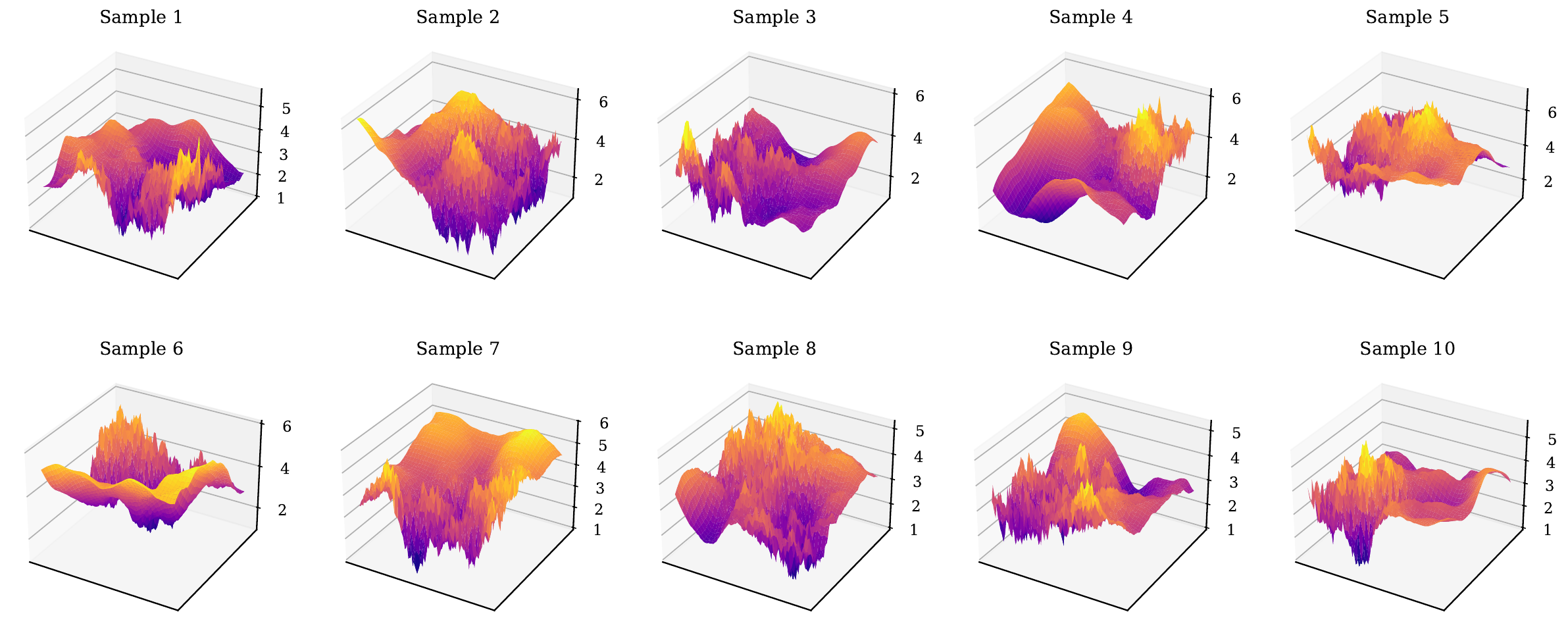}
\caption{Batch of $10$ inhomogeneous benchmark functions with randomly generated anchor points $\boldsymbol{\mu}_i$.}
\label{fig:GP_grid}
\end{figure}
To validate the consistency of these improvements, we generated a benchmark set of $10$ inhomogeneous functions. While the kernel types remained constant, the three anchor points $\boldsymbol{\mu}_i$ determining the spatial zones for each function were randomly generated with a minimum distance constraint of $0.3$ between any pair of points to ensure diverse spatial configurations (see Figure~\ref{fig:GP_grid}).

\begin{figure}[t]
\centering
\includegraphics[width=\textwidth]{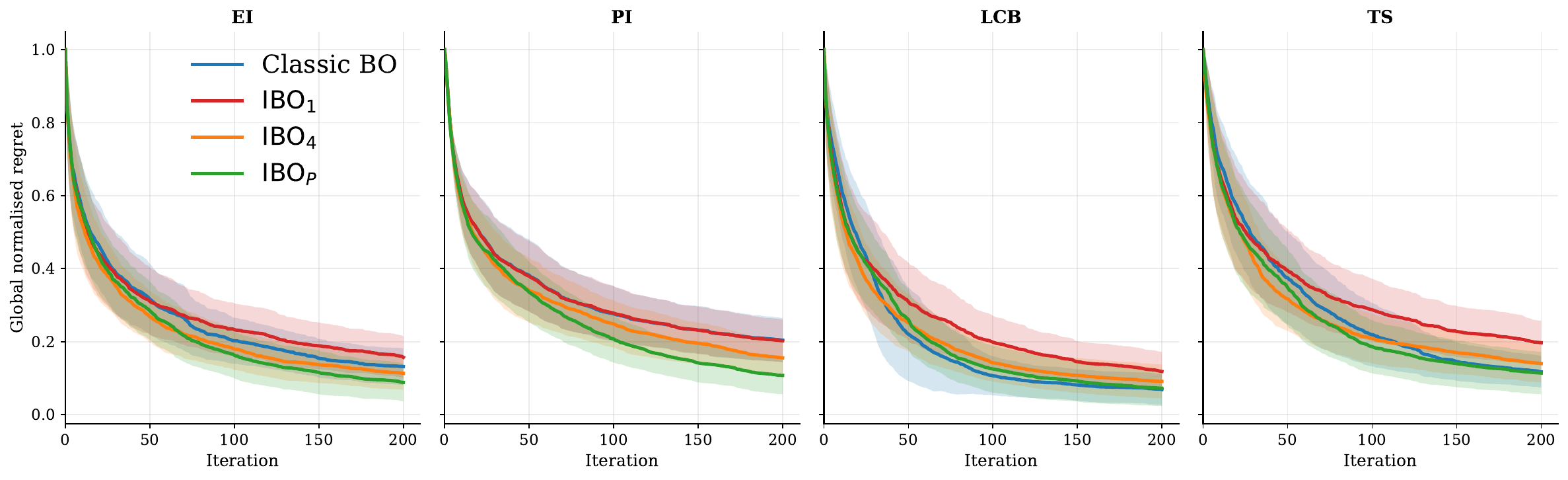}
\caption{Global normalized regret on the $10$ inhomogeneous functions averaged over the functions. Shaded regions are $95\%$ confidence intervals across functions.}
\label{fig:nstat_ratios}
\end{figure}

Figure~\ref{fig:nstat_ratios} reports the global normalized regret~\eqref{eq:global_nr} on these $10$ functions. $\text{IBO}_P$ outperforms Classic BO under EI and PI, while the four methods are more or less indistinguishable under LCB and TS. The criterion partitions every surface in almost every run, with $5.17$ leaves on average, and the mean differences of Table~\ref{tab:main} confirm that the advantage under EI and PI is significant.

Table~\ref{tab:ablation} isolates where the gain comes from. Estimating the scale is not enough on its own: $\mathrm{IBO}_1$, which estimates a single scale on the whole domain, is beaten by Classic BO on $31$ of the $40$ function--acquisition pairs, with a median ratio of final normalized regrets of $0.81$. Partitioning reverses this, and the adaptive partition multiplies the gain: $\mathrm{IBO}_4$ reaches a median ratio of $1.08$ and $\mathrm{IBO}_P$ of $1.53$.

\begin{table}[t]
\centering
\caption{Partitioning statistics of $\mathrm{IBO}_P$ and comparison with Classic BO. Each acquisition row gives the mean over the $300$ function--seed pairs of the difference in final normalized regret between the two methods at the same seed, $\mathrm{IBO}_P$ minus Classic BO (negative values favor $\mathrm{IBO}_P$), with a star when its $95\%$ confidence interval excludes zero, followed by the number of pairs on which $\mathrm{IBO}_P$\,/\,Classic BO attains the smaller final normalized regret (ties are not counted).}

\label{tab:main}
\begin{tabular}{lcc}
\toprule
& 10 homogeneous functions & 10 inhomogeneous functions \\
\midrule
Mean number of leaves         & $4.36$  & $5.17$   \\
Runs with a split             & $68\%$  & $99.8\%$ \\
\midrule
EI  & $+0.018$\phantom{$^*$}\ \ (124/133) & $-0.043^*$\ \ (209/90)  \\
PI  & $-0.041^*$\ \ (157/47)              & $-0.097^*$\ \ (227/70)  \\
LCB & $+0.056^*$\ \ (121/145)             & $+0.004$\phantom{$^*$}\ \ (132/168) \\
TS  & $+0.011$\phantom{$^*$}\ \ (143/90)  & $-0.003$\phantom{$^*$}\ \ (193/107) \\
\bottomrule
\end{tabular}
\end{table}

\begin{table}[t]
\centering
\caption{Comparison of IBO against Classic BO: number of functions, out of $10$, on which each IBO variant attains the smaller normalized regret~\eqref{eq:nr} at the final step, and the median ratio of the final normalized regret of Classic BO to that of the variant. A median ratio above one means that IBO
outperforms BO.}
\label{tab:ablation}
\begin{tabular}{lccc@{\hskip 2em}ccc}
\toprule
 & \multicolumn{3}{c}{10 homogeneous functions}
 & \multicolumn{3}{c}{10 inhomogeneous functions} \\
\cmidrule(lr){2-4}\cmidrule(lr){5-7}
 & $\mathrm{IBO}_1$ & $\mathrm{IBO}_4$ & $\mathrm{IBO}_P$
 & $\mathrm{IBO}_1$ & $\mathrm{IBO}_4$ & $\mathrm{IBO}_P$ \\
\midrule
EI    & 0 & 2 & \textbf{5} & 1 & \textbf{8} & \textbf{8} \\
PI    & 3 & 5 & \textbf{8} & 5 & 8 & \textbf{9} \\
LCB   & 1 & 3 & \textbf{5} & 2 & 4 & \textbf{5} \\
TS    & 1 & 3 & \textbf{5} & 1 & 3 & \textbf{5} \\
\midrule
Total & $5/40$ & $13/40$ & $\mathbf{23/40}$
      & $9/40$ & $23/40$ & $\mathbf{27/40}$ \\
Median ratio & $0.71$ & $0.76$ & $\mathbf{1.12}$
             & $0.81$ & $1.08$ & $\mathbf{1.53}$ \\
\bottomrule
\end{tabular}
\end{table}

\subsection{Real-World Application: LASSO 8-Dimensional Optimization}

To assess the practical utility of the proposed framework, we consider the hyperparameter optimization of a LASSO model applied to the diabetes dataset from the LassoBench suite~\cite{vsehic2022lassobench}, an $8$-dimensional
search space that serves as a standard benchmark for sparse regression. The task involves tuning the regularization path and associated model parameters to minimize the validation error, over a budget of $200$ iterations.

\begin{figure}[t]
\centering
\includegraphics[width=\textwidth]{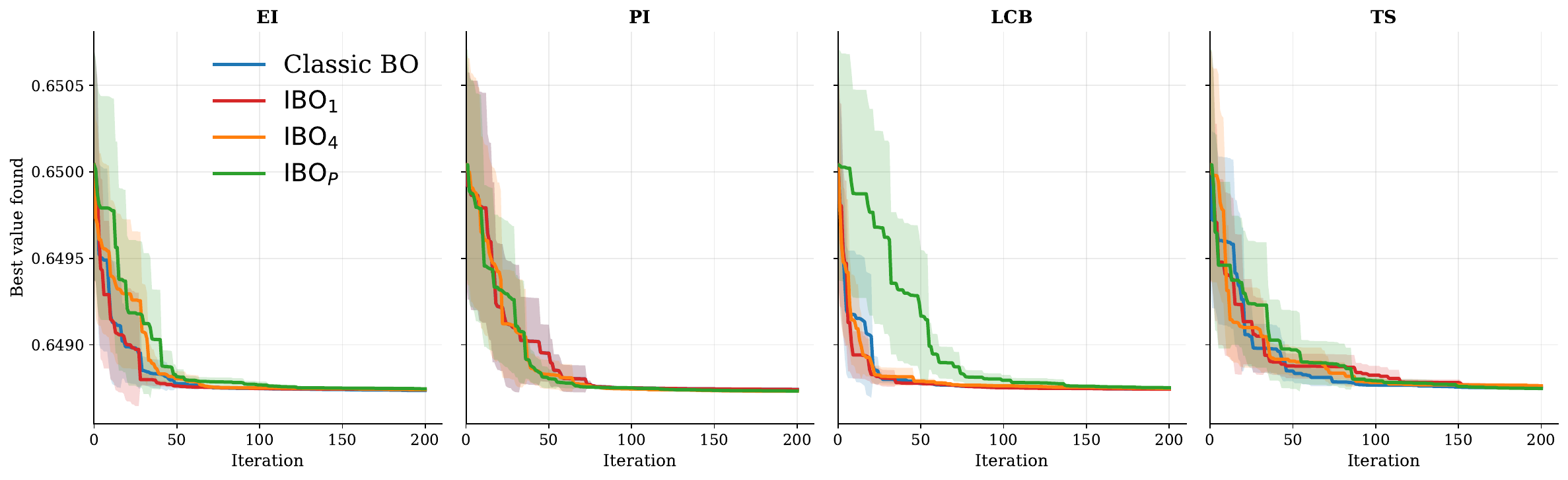}
\caption{Mean convergence curves on the LASSO task. Shaded areas represent the 95\% confidence interval across seeds.}
\label{fig:lasso_results}
\end{figure}

Convergence results are summarized in Figure~\ref{fig:lasso_results}. The objective varies by less than $0.3\%$ over the whole run, from $0.6505$ at the initial design to $0.6487$ at the end, and all four methods attain that value
within the budget. Under EI and PI they are indistinguishable after fifty evaluations; under LCB, $\text{IBO}_P$ needs about seventy more evaluations than the others to reach the same value, and under TS all four
have converged by roughly a hundred.

On this task the four methods reach closely comparable final objective values, and we draw a correspondingly narrow conclusion: partitioning does not degrade final optimization performance. Whether this is a case in which
the framework should help is not something we can determine, since the regularity of the LassoBench validation surface is unknown to us, in contrast to the synthetic functions of Section~\ref{sec:nonstat}, which are constructed to be inhomogeneous.

\section{Conclusions}
\label{sec:conclusions}

We introduced the IBO framework to address the limitations of stationary Gaussian processes when optimizing objective functions with spatially varying regularity. By moving beyond the assumption of a spatially homogeneous latent function, we developed a partitioned modeling approach that adaptively adjusts the spatial resolution of the surrogate model using locally adaptive kernels.
Our theoretical analysis established that polynomially denser sampling in rougher subdomains suffices to certify regional acquisition suprema. The primary algorithmic contribution of this work was the recursive splitting criterion in Section~\ref{subsec:RPHAO}. Driven by the non-stationarity ratio $\eta$, this criterion effectively learns the spatially varying effective roughness of the objective function. Empirical results across inhomogeneous functions and an $8$-dimensional real-world task show that the IBO framework reduces regret where the local regularity varies, under the improvement-based acquisitions, while achieving mixed results on homogeneous functions.
While the IBO framework provides a principled strategy for navigating inhomogeneous objective functions, possibilities remain for future extensions. In particular, integrating this framework with manifold learning~\cite{borovitskiy2020matern, jaquier2020high} would allow IBO to handle inhomogeneous functions defined on non-Euclidean search spaces, providing a grounded, adaptive strategy for navigating complex, heterogeneous surfaces.

\section*{Acknowledgments}
TK acknowledges the research environment provided by ELLIS Institute Finland. The authors used Claude for proofreading, polishing the text, and for code review. The authors assume responsibility for all content.

\bibliographystyle{siamplain}
\bibliography{references}

\end{document}